\documentclass[11pt]{amsart}
\usepackage{amssymb,amsmath,amsgen,amsxtra,amsfonts} 
\usepackage{amscd}
\usepackage{dsfont,times}

\begin{document}
%
%
%
\theoremstyle{definition}
\newtheorem{Definition}{Definition}[section]
\newtheorem*{Definitionx}{Definition}
\newtheorem{Convention}{Definition}[section]
\newtheorem{Construction}{Construction}[section]
\newtheorem{Example}[Definition]{Example}
\newtheorem{Examples}[Definition]{Examples}
\newtheorem{Remark}[Definition]{Remark}
\newtheorem{Remarks}[Definition]{Remarks}
\newtheorem{Caution}[Definition]{Caution}
\newtheorem{Conjecture}[Definition]{Conjecture}
\newtheorem*{Conjecturex}{Conjecture}
\newtheorem{Question}[Definition]{Question}
\newtheorem{Questions}[Definition]{Questions}
\newtheorem*{Acknowledgements}{Acknowledgements}
\theoremstyle{plain}
\newtheorem{Theorem}[Definition]{Theorem}
\newtheorem*{Theoremx}{Theorem}
\newtheorem{Proposition}[Definition]{Proposition}
\newtheorem*{Propositionx}{Proposition}
\newtheorem{Lemma}[Definition]{Lemma}
\newtheorem{Corollary}[Definition]{Corollary}
\newtheorem*{Corollaryx}{Corollary}
\newtheorem{Fact}[Definition]{Fact}
\newtheorem{Facts}[Definition]{Facts}
\newtheoremstyle{voiditstyle}{3pt}{3pt}{\itshape}{\parindent}%
{\bfseries}{.}{ }{\thmnote{#3}}%
\theoremstyle{voiditstyle}
\newtheorem*{VoidItalic}{}
\newtheoremstyle{voidromstyle}{3pt}{3pt}{\rm}{\parindent}%
{\bfseries}{.}{ }{\thmnote{#3}}%
\theoremstyle{voidromstyle}
\newtheorem*{VoidRoman}{}

%
\newcommand{\prf}{\par\noindent{\sc Proof.}\quad}
\newcommand{\blowup}{\rule[-3mm]{0mm}{0mm}}
\newcommand{\cal}{\mathcal}
\newcommand{\Aff}{{\mathds{A}}}
\newcommand{\BB}{{\mathds{B}}}
\newcommand{\CC}{{\mathds{C}}}
\newcommand{\FF}{{\mathds{F}}}
\newcommand{\GG}{{\mathds{G}}}
\newcommand{\HH}{{\mathds{H}}}
\newcommand{\NN}{{\mathds{N}}}
\newcommand{\ZZ}{{\mathds{Z}}}
\newcommand{\PP}{{\mathds{P}}}
\newcommand{\QQ}{{\mathds{Q}}}
\newcommand{\RR}{{\mathds{R}}}
\newcommand{\Liea}{{\mathfrak a}}
\newcommand{\Lieb}{{\mathfrak b}}
\newcommand{\Lieg}{{\mathfrak g}}
\newcommand{\Liem}{{\mathfrak m}}
\newcommand{\ideala}{{\mathfrak a}}
\newcommand{\idealb}{{\mathfrak b}}
\newcommand{\idealg}{{\mathfrak g}}
\newcommand{\idealm}{{\mathfrak m}}
\newcommand{\idealp}{{\mathfrak p}}
\newcommand{\idealq}{{\mathfrak q}}
\newcommand{\idealI}{{\cal I}}
\newcommand{\lin}{\sim}
\newcommand{\num}{\equiv}
\newcommand{\dual}{\ast}
\newcommand{\iso}{\cong}
\newcommand{\homeo}{\approx}
\newcommand{\mm}{{\mathfrak m}}
\newcommand{\pp}{{\mathfrak p}}
\newcommand{\qq}{{\mathfrak q}}
\newcommand{\rr}{{\mathfrak r}}
\newcommand{\pP}{{\mathfrak P}}
\newcommand{\qQ}{{\mathfrak Q}}
\newcommand{\rR}{{\mathfrak R}}
%
%
\newcommand{\dq}{{``}}
\newcommand{\OO}{{\cal O}}
\newcommand{\into}{{\hookrightarrow}}
\newcommand{\onto}{{\twoheadrightarrow}}
\newcommand{\Spec}{{\rm Spec}\:}
\newcommand{\BigSpec}{{\rm\bf Spec}\:}
\newcommand{\Proj}{{\rm Proj}\:}
\newcommand{\Pic}{{\rm Pic }}
\newcommand{\Br}{{\rm Br}}
\newcommand{\NS}{{\rm NS}}
\newcommand{\chit}{\chi_{\rm top}}
\newcommand{\KanDiv}{{\cal K}}
\newcommand{\Cycl}[1]{{\ZZ/{#1}\ZZ}}
\newcommand{\Sym}{{\mathfrak S}}
\newcommand{\ab}{{\rm ab}}
\newcommand{\Aut}{{\rm Aut}}
\newcommand{\Hom}{{\rm Hom}}
\newcommand{\ord}{{\rm ord}}
\newcommand{\divisor}{{\rm div}}
\newcommand{\Alb}{{\rm Alb}}
\newcommand{\Jac}{{\rm Jac}}
\newcommand{\Km}{{\rm Km}}
\newcommand{\Ig}[1]{{\rm Ig}(#1)}
\newcommand{\Igord}[1]{{\rm Ig}(#1)^{\rm ord}}
\newcommand{\oIgord}[1]{{\overline{{\rm Ig}(#1)^{\rm ord}}}}
\newcommand{\Brauer}{\widehat{{\rm Br}}}
\newcommand{\MW}{{\rm MW}}
\newcommand{\nMW}{{\rm MW}^\circ}
\newcommand{\piet}{{\pi_1^{\rm \acute{e}t}}}
\newcommand{\Het}[1]{{H_{\rm \acute{e}t}^{{#1}}}}
\newcommand{\Hcris}[1]{{H_{\rm cris}^{{#1}}}}
\newcommand{\HdR}[1]{{H_{\rm dR}^{{#1}}}}
\newcommand{\hdR}[1]{{h_{\rm dR}^{{#1}}}}
\newcommand{\defin}[1]{{\bf #1}}

\title[Elliptic K3 surfaces]{Elliptic K3 surfaces with ${\mathbf p}^{\mathbf n}$-torsion sections}


\author[Hiroyuki Ito]{Hiroyuki Ito}
\address{Department of Applied Mathematics,
Graduate School of Engineering, Hiroshima University, Higashi-Hiroshima 739-8527, Japan}
\curraddr{}
\email{hiroito@amath.hiroshima-u.ac.jp}

\author[Christian Liedtke]{Christian Liedtke}
\address{Department of Mathematics, Stanford University, 450 Serra Mall, Stanford CA 94305-2125, USA}
\curraddr{}
\email{liedtke@math.stanford.edu}

\dedicatory{February 4, 2011}
\maketitle

\begin{abstract}
 We classify elliptic K3 surfaces in characteristic $p$ with
 $p^n$-torsion section.
 For $p^n\geq3$ we verify conjectures of Artin and Shioda, compute
 the heights of their formal Brauer groups, as well as Artin invariants 
 and Mordell--Weil groups in the supersingular cases.
\end{abstract}

\section*{Introduction}

The geometry and arithmetic of K3 surfaces is a fascinating subject of
algebraic geometry.
Moreover, this class of surfaces provides a rich source of conjectures 
that are difficult to come by.

In this paper, we consider K3 surfaces in positive characteristic $p$
that are elliptically fibered.
Moreover, we assume that the fibration possesses a torsion section
of order $p^n$.
Such surfaces have already been studied by Schweizer \cite{Schweizer 2005}.
Recall, e.g., from \cite[Chapter 12]{Katz; Mazur 1985},
that the Igusa moduli functor, which classifies ordinary
elliptic curves with $p^n$-torsion sections, is representable by
a smooth affine curve, the so-called Igusa curve $\Igord{p^n}$
if $p^n\geq3$.
Using Igusa's results \cite{Igusa 1968},
we first strengthen results of \cite{Schweizer 2005} and
\cite{Dolgachev Keum 2009}:

\begin{Theoremx}
  Elliptic K3 surfaces with $p^n$-torsion section in characteristic 
  $p$ exist for $p^n\leq8$ only.
  If the fibration has constant $j$-invariant then $p^n=2$.
\end{Theoremx}

Using the universal elliptic curves over the Igusa curves and
the results \cite{Liedtke; Schroeer 2008}
on their N\'eron models over their cusps and the
supersingular locus, we explicitly classify elliptic K3
surfaces with $p^n$-torsion sections for $p^n\geq3$.

Next, translation by a $p$-torsion section of an elliptic 
fibration induces a 
$\ZZ/p\ZZ$-action, i.e., a wild $p$-cyclic automorphism.
Such wild automorphisms on K3 surfaces have been studied 
in general by Dolgachev and Keum \cite{Dolgachev; Keum 2001}.
Using their results, we illustrate and strengthen these
results in case the wild automorphism arises from translation
by a $p$-torsion section.
For example, we determine the fixed point set of translation
by a $p$-torsion section in bad fibers of the elliptic fibration,
which extends work of Miranda and Persson \cite{Miranda; Persson 1989}
from the prime-to-$p$ case.

Before stating one of our main results, let us state
a couple of conjectures on the arithmetic of elliptic K3 surfaces.
First, let us recall that a surface is called
{\em Shioda-supersingular} if the rank of its N\'eron--Severi
group is equal to its second Betti number.
In \cite{Shioda 1974b}, Shioda has shown that unirational surfaces are 
Shioda-supersingular and conjectured the converse
in \cite{Shioda 1977}.
On the other hand, a surface is called
{\em Artin-supersingular} if its formal Brauer group has 
infinite height.
Artin \cite{Artin 1974} has shown that unirational K3 surfaces
are Artin-supersingular and conjectured the converse.
Moreover, he proved in loc. cit. that Shioda-supersingular
surfaces are Artin-supersingular and conjectured the
converse.
Thus

\begin{Conjecturex}
 For K3 surfaces,
 \begin{enumerate}
  \item (Shioda) Shioda-supersingularity implies unirationality,
  \item (Artin) Artin-supersingularity implies unirationality,
  \item (Artin) Artin-supersingularity implies Shioda-supersingularity.
 \end{enumerate}
\end{Conjecturex}

For elliptic K3 surfaces these two notions of supersingularity
coincide \cite{Artin 1974}.
In characteristic $2$,
there is another conjecture by Artin
\cite{Artin 1974}, 
which does not only imply the above conjectures but also
gives a geometric explanation of the above conjectures:

\begin{Conjecturex}[Artin]
 In characteristic $2$, an elliptic fibration on a supersingular 
 K3 surface
 arises via Frobenius pullback from a rational elliptic surface.
\end{Conjecturex}

Unfortunately, such a conjecture cannot be true in general
in characteristic $p\geq3$, see Section \ref{sec:brauer} for discussion.
However, for elliptic K3 surfaces with $p^n$-torsion sections
a beautiful picture emerges:

\begin{Theoremx}
 Let $X\to\PP^1$ be an elliptic K3 surface 
 in positive characteristic $p$ with $p^n$-torsion 
 sections and $p^n\geq3$.
 Then the following are equivalent
 \begin{enumerate}
  \item The elliptic fibration arises as Frobenius pullback from
    a rational elliptic fibration.
  \item $X$ is unirational.
  \item $X$ is supersingular.
  \item The fibration has precisely one additive fiber.
 \end{enumerate}
\end{Theoremx}

\begin{Corollaryx}
 The conjectures of Artin and Shioda hold for
 elliptic K3 surfaces with $p^n$-torsion sections if
 $p^n\geq3$.
\end{Corollaryx}

Let us recall that the moduli space of K3 
surfaces is stratified by the height $h$ of the formal 
Brauer group, which takes every value $1\leq h\leq 10$ or 
$h=\infty$.
Furthermore, the moduli space of surfaces with 
$h=\infty$, i.e., the Artin-supersingular surfaces, is stratified
by the Artin invariant $\sigma_0$, which takes
every value $1\leq\sigma_0\leq10$.
For our surfaces we prove the following alternative

\begin{Propositionx}
 For an elliptic K3 surface with $p^n$-torsion section
 in characteristic $p$ and $p^n\geq3$ there are two
 possibilities:  
 \begin{enumerate}
  \item either the elliptic fibration has precisely one additive
   fiber and the surface is supersingular ($h=\infty$), 
  \item or the elliptic fibration has precisely two additive fibers
   and the surface is ordinary ($h=1$)
 \end{enumerate}
\end{Propositionx}

In characteristic $2$, a connection between the height of the
formal Brauer group and the singular fibers of an elliptic
fibration has already been observed by Artin \cite{Artin 1974}.
For $p^n\geq7$ there is only one elliptic K3 surface
with $p^n$-torsion section and it is supersingular.
On the other hand, the generic elliptic K3 surface with
$p^n$-torsion section with $p^n\leq5$ is ordinary.

Concerning the Artin invariants of the supersingular surfaces
we obtain the following characterization:

\begin{Theoremx}
  The Artin invariant $\sigma_0$ of a supersingular and elliptic K3 surface with
  $p^n$-torsion in characteristic $p$ satisfies $\sigma_0\leq\sigma_0(p^n)$ where
  $$\begin{array}{c|ccccc}
    p^n & 8 & 7 & 5 & 4 & 3\\
    \hline
    \sigma_0(p^n) & 1 & 1 & 2 & 3 & 6
   \end{array}$$
 Conversely, a supersingular K3 surface in characteristic $p$
 with $\sigma_0\leq\sigma_0(p^n)$ possesses an elliptic fibration
 with $p^n$-torsion section.
\end{Theoremx}

We also determine the Mordell--Weil groups and find 
explicit Weierstra\ss\ equations of these fibrations.
In particular, we obtain explicit and complete families
of supersingular K3 surfaces with $\sigma_0\leq\sigma_0(p^n)$
in characteristic $p$.
To obtain these results in characteristic $p\leq3$ 
we use semi-universal deformations of the $E_8^2$-singularity 
($p=3$) and the $E_8^4$-singularity ($p=2$).
\medskip
%

On the other hand, elliptic K3 surfaces with $2$-torsion section
in characteristic $2$ are much harder to come by.
This has to do with the fact that there is no Igusa curve
to ''tame`` the situation.
It turns out that there are extra classes.
For example, fibrations with constant $j$-invariant have to be considered
and there are classes where the formal Brauer group has
height $2$, i.e., the above alternative does no longer hold.
We refer to Theorem \ref{p=2theorem} for the precise 
structure result.

%
\medskip

The article is organized as follows:
In Section \ref{sec:generalities} we recall a couple of general facts
about the Igusa moduli problem and show that elliptic K3 surfaces
with $p^n$-torsion sections can exist for $p^n\leq8$ only.
In Section \ref{sec:wild} we analyze the fixed locus of translation
by a $p$-torsion section in an elliptic fibration.
In Section \ref{sec:brauer} we compute the height of the formal
Brauer group in terms of the additive fibers of an elliptic
fibration.
This already yields some of our main theorems for $p\neq2$.
In Section \ref{sec:classification} we give an explicit classification
for $p\neq2$ and compute the Artin invariants in the supersingular
cases.
The rest of the article takes place in characteristic $2$ only:
in Section \ref{sec:p=2} we prove the general structure result, and
classify the new, ''exotic`` classes in Section \ref{sec:newp=2}.
Finally, in Section \ref{sec:p=4} we deal with $4$- and $8$-torsion sections
and use again the corresponding Igusa curves.

\begin{Acknowledgements}
 We thank Igor~Dolgachev, Matthias~Sch\"utt and the referee 
 for helpful comments.
 The first author acknowledges the support by Grant-in-Aid for
 Scientific Research (C) 20540044, the Ministry of Education, Culture,
 Sports, Science and Technology.
 The second author gratefully acknowledges funding from DFG under 
 research grant LI 1906/1-1 and thanks the department of mathematics
 at Stanford university for kind hospitality.
\end{Acknowledgements}

\section{Igusa curves}
\label{sec:generalities}

In this section we first recall the Igusa moduli problem and the
Igusa curves and use these results to show that 
elliptic K3 surfaces with $p^n$-torsion section can exist only 
if $p^n\leq8$.

Let us recall, e.g. from \cite[Chapter 12.3]{Katz; Mazur 1985}, 
that the Igusa moduli functor $[\Igord{p^n}]$ 
associates to every scheme $S$ over $\FF_p$ 
the set of ordinary elliptic curves 
${\cal E}$ over $S$ such that the $n$-fold Frobenius pullback
${\cal E}^{(p^n)}=(F^n)^*({\cal E})$ possesses a $p^n$-torsion section.
If $p^n\geq3$ then this functor is representable by a 
smooth and affine curve over $\FF_p$, the {\em Igusa curve}
$\Igord{p^n}$.
We denote by ${\cal E}\to\Igord{p^n}$ the universal family.
Thus, if $X\to B$ is an elliptic fibration in characteristic
$p$ with $p^n$-torsion section, and if $U\subseteq B$ 
denotes the open set over which the fibres are ordinary
elliptic curves, then
there exists a classifying morphism
$$
\varphi\,:\,U\,\to\,\Igord{p^n}
$$
such that the restriction $X|_U\to U$ is isomorphic
to $(F^n)^*({\cal E})\to U$.

The geometry of the normal compactification $\oIgord{p^n}$ of
$\Igord{p^n}$ has been studied in \cite{Igusa 1968}.
For example, if $n=1$ and $p\geq3$, which is the
case that we will be needing most in the sequel,
then the $j$-invariant induces a Galois morphism
$\oIgord{p^n}\to\PP^1$, whose Galois group is cyclic 
of order $(p-1)/2$.
This morphism is totally ramified over the supersingular 
$j$-values and totally split over $j=\infty$, i.e.,
there are $(p-1)/2$ points lying above infinity,
the so-called cusps.
The degenerating behavior of the universal family
${\cal E}\to\Igord{p}$ over the supersingular points
and the cusps has been determined
in \cite{Liedtke; Schroeer 2008}.

\begin{Theorem}
 \label{possible cases}
 An elliptic K3 surface $X\to\PP^1$ with $p^n$-torsion sections in
 positive characteristic $p$ satisfies the inequality $p^n\leq8$.
 Moreover, if the fibration has constant $j$-invariant then $p=2$ and $n=1$.
\end{Theorem}

\prf
We first deal with the case of constant $j$-invariant.
Since the $p^n$-torsion section is different from the
zero section, the generic fiber is ordinary and so
the ordinary locus $U\subseteq\PP^1$ is open and dense.
Moreover, if $p^n\geq3$ then the Igusa moduli problem is representable
and constant $j$-invariant implies that the classifying morphism
$\varphi:U\to\Igord{p^n}$ is constant.
Thus, $X|_U\to U$ is a product family
(the Igusa curve is a fine moduli space),
and not birational to a K3 surface.
Hence in this case we have $p^n=2$.

We may thus assume that the fibration has non-constant
$j$-invariant, and again, 
the ordinary locus $U\subseteq\PP^1$ is open and dense.
Also, we may assume $p^n\geq3$, i.e., that the Igusa moduli problem
is representable.
Then the classifying morphism $\varphi:U\to\Igord{p^n}$
is dominant, which implies that
$\Igord{p^n}$ is a rational curve.
The genera of the Igusa curves have been determined in
\cite{Igusa 1968} and a straightforward computation shows
that these curves are rational if and only if $p^n\leq11$.

Let us first exclude $p=11$.
In this case $\oIgord{11}$ has $5$ cusps.
Hence our fibration has at least $5$ fibres with
potentially multiplicative reduction.
By \cite[Theorem 4.3]{Liedtke; Schroeer 2008}
we have in fact multiplicative reduction.
Thus, our family
has at least $5$ fibres with multiplicative 
reduction, necessarily of type ${\rm I}_n$, where $11$
divides all these $n$'s.
These contribute at least $5\times(p-1)=50$ to $\rho(X)$,
i.e., $b_2(X)\geq\rho(X)>50$, i.e., $X$ is not a K3 surface.

The remaining case $p^n=9$ is excluded similarly and
we leave it to the reader.
\qed\medskip

\begin{Remark}
  Non-existence of elliptic K3 surfaces with $p$-torsion sections
  for $p\geq11$ has been shown in  
  \cite[Theorem 2.13]{Dolgachev Keum 2009}.
  Under the assumption that the fibration does not have
  constant $j$-invariant, Theorem \ref{possible cases} 
  has been shown in the remark after 
  \cite[Theorem 2.3]{Schweizer 2005}, using methods closely
  related to ours.
  Nevertheless, we decided to give a proof in our setup, i.e.,
  by analyzing the classifying morphisms to the Igusa curves
  and their universal families.
\end{Remark}

The proof shows that Igusa curves that are rational are 
crucial for the description
of elliptic K3 surfaces with $p^n$-torsion sections.
Igusa's results \cite{Igusa 1968} show that these
curves are rational if and only if $p^n\leq11$.
For our explicit classification later on, and
in order to obtain equations when needed, we 
determine Weierstra\ss\ equations in these cases.

\begin{Proposition}
 \label{explicit equations}
 The universal elliptic curves over $\Igord{p^n}$ for
 $p^n\leq11$ are given by the following 
 equations over $\FF_p[t]$:
 $$\begin{array}{clclcc}
    \blowup
    p^n & & & &\multicolumn{2}{c}{\mbox{ singular fibres }}\\
    \hline
    11 & {\cal E} &:& y^2 \,=\,x^3\,+\,(t-1)^{-1}t x\,+\,5t^{-1}(t-1) & 5\times{\rm I}_1, & {\rm II}^*, {\rm III}^* \\
       & {\cal E}^{(p)} &:& y^2 \,=\,x^3\,+\,(t-1)^{-11} t^{11} x\,+\,5t^{-11} (t-1)^{11} 
                        & 5\times{\rm I}_{11}, & {\rm II}, {\rm III} \\
    9  & {\cal E}         &:& y^2+txy=x^3-t^3(t^2-1)& 3\times{\rm I}_1, & {\rm IV}^*_1 \\
       & {\cal E}^{(p)} &:& y^2+tx+(t^2-t)y=x^3+tx+(t^2-t) & 3\times{\rm I}_3, & {\rm II}_1 \\
       & {\cal E}^{(p^2)} &:& y^2+t^3x+(t^6-t^3)y=x^3+t^3x+(t^6-t^3) & 3\times{\rm I}_9, & {\rm IV}^*_1\\
    8  & {\cal E}         &:& y^2+xy=x^3+t(t+1)     & 2\times{\rm I}_1, & {\rm III}^*_1\\
       & {\cal E}^{(p)}   &:& y^2+xy=x^3+t^2(1+t^2) & 2\times{\rm I}_2, & {\rm I}^*_{1,1}\\
       & {\cal E}^{(p^2)} &:& y^2+xy=x^3+t^4(1+t^4) & 2\times{\rm I}_4, & {\rm III}_1\\
       & {\cal E}^{(p^3)} &:& y^2+xy=x^3+t^8(1+t^8) & 2\times{\rm I}_8, & {\rm I}^*_{1,1}\\
    7  & {\cal E}       &:& y^2=x^3+t^3x+5t^6    & 3\times{\rm I}_1, & {\rm III}^* \\
       & {\cal E}^{(p)} &:& y^2=x^3+tx+5t^{12}   & 3\times{\rm I}_7, & {\rm III}\\
    5  & {\cal E}       &:& y^2=x^3+3t^4x+t^5    & 2\times{\rm I}_1, & {\rm II}^*\\
       & {\cal E}^{(p)} &:& y^2=x^3+3t^4x+t      & 2\times{\rm I}_5, & {\rm II}\\
    4  & {\cal E}         &:& y^2+xy=x^3+t       & {\rm I}_1, & {\rm II}^*_1 \\
       & {\cal E}^{(p)}   &:& y^2+xy=x^3+t^2     & {\rm I}_2, & {\rm III}^*_1 \\
       & {\cal E}^{(p^2)} &:& y^2+xy=x^3+t^4     & {\rm I}_4, & {\rm I}^*_{1,1}\\
    3  & {\cal E}       &:& y^2+txy=x^3-t^5      & {\rm I}_1, & {\rm II}^*_1 \\
       & {\cal E}^{(p)} &:& y^2+txy+t^2y=x^3     & {\rm I}_3, & {\rm IV}^*_1 \\
   \end{array}
 $$
 All places of bad reduction are defined over $\FF_p$ with split multiplicative
 reduction at the cusps and additive reduction at the supersingular points.
\end{Proposition}

\proof
As an example we do the case $p=7$ and leave the others to the reader:
The elliptic curve ${\cal E}$ for $p=7$ given in the table has Hasse invariant 
$[1]\in\FF_p/\FF_p^{\times(p-1)}$, which implies that ${\cal E}^{(p)}$ has a
$\FF_p[t]$-rational $p$-division point.
Thus, there exists a morphism $\varphi:\Spec\FF_p[t]\to\oIgord{p}$ such that
$\cal E$ is the pullback of the universal elliptic curve over $\Igord{p}$
via $\varphi$.
Since the $j$-invariant of $\cal E$ is not constant, it follows that 
$\varphi$ is a finite morphism.
The curve $\oIgord{p}$ has $(p-1)/2=3$ cusps over which the
universal family degenerates into ${\rm I}_1$-fibers
\cite[Theorem 10.3]{Liedtke; Schroeer 2008}.
Since the same is true for ${\cal E}$, we get $\deg\varphi=1$, i.e.,
$\varphi$ is an isomorphism.
\qed

\section{Wild $p$-cyclic actions}
\label{sec:wild}

Since we are dealing with elliptic fibrations with
$p$-torsion sections in positive characteristic $p$, 
translation by such a torsion section gives rise to a wild automorphism,
and we may apply the results of \cite{Dolgachev; Keum 2001}.
For K3 surfaces, we will see that there are at most two additive fibers and 
if there are two such fibers then the elliptic fibration 
arises as Frobenius pullback from an elliptic K3 surface.
To fix notation, let $X\to B$ be an elliptic surface 
with zero section $\sigma_0$ and $p$-torsion section $\sigma_p$.
We denote by $G$ the cyclic group of order $p$ 
generated by translations by $\sigma_p$ and set $Y:=X/G$.
Note that the elliptic fibration $X\to B$ induces
an elliptic fibration $Y\to B$ and we get a diagram
of elliptic fibrations over $B$
\begin{equation}
 \label{multiplication by p}
 Y \,\to\, X \,\to\, X/G\,\iso\,Y\,,
\end{equation}
where the first map is purely inseparable (relative Frobenius over $B$)
and the second is an Artin-Schreier morphism.

We now analyze the action of $G$ induced on the fibers.
In characteristic zero and for multiplicative reduction, this has been 
worked out in \cite[Section 2]{Miranda; Persson 1989}.
If $X_0$ denotes a special fiber of the fibration we will denote
by $(\sigma_0\cdot\sigma_p)_0$ the intersection number of $\sigma_0$
and $\sigma_p$ in the fiber $X_0$.
Finally, we denote by $F_0$ the reduced fixed point scheme
of the $G$-action on $X_0$, see also the discussion in
\cite[Remark 2.7]{Dolgachev; Keum 2001}.

\begin{Proposition}
  \label{fixed locus of translation}
  Let $X\to B$ be an elliptic fibration in characteristic $p$ with
  $p$-torsion section $\sigma_p$.
  Let $X_0$ be a special fiber and let $F_0$ be the reduced fixed point 
  scheme of the $\sigma_p$-translation on $X_0$.

  If $X_0$ has semi-stable reduction and more precisely, if the reduction is
  \begin{enumerate}
    \item good and ordinary then $(\sigma_0\cdot\sigma_p)_0=0$
       and $F_0=\emptyset$,
    \item good and supersingular then $(\sigma_0\cdot\sigma_p)_0\geq1$ 
       and $F_0=X_0$,
    \item bad multiplicative then $(\sigma_0\cdot\sigma_p)_0=0$ 
       and $F_0=\emptyset$.
  \end{enumerate}

  If $X_0$ has additive reduction and $(\sigma_0\cdot\sigma_p)_0\geq1$ then
  $F_0=X_0$.

  If $X_0$ has additive reduction, $(\sigma_0\cdot\sigma_p)_0=0$ and 
  the reduction type is
  \begin{enumerate}
    \item ${\rm II}$, ${\rm III}$, ${\rm IV}$  then $F_0$ equals the unique 
       point that is not smooth over the base of the fibration,
    \item ${\rm I}_n^*$ ($p\neq2$), ${\rm IV}^*$ ($p\neq3$), 
       ${\rm III}^*$ ($p\neq2$), ${\rm II}^*$
       then $F_0$ is a curve, equal to the union of all multiplicity 
       $\geq2$-components of $X_0$
  \end{enumerate}
  In characteristic $p\leq3$ the situation is the same if $\sigma_p$ does not 
  specialize into the component group of $X_0$. 
  If it does and if the reduction type is
  \begin{enumerate}
     \item ${\rm IV}^*$ then $p=3$ and $F_0$ is one point, which lies on the 
       component of multiplicity $3$,
     \item ${\rm III}^*$ then $p=2$ and $F_0$ is one point, namely the intersection
       of the component of multiplicity $4$ and the one of multiplicity $2$,
     \item ${\rm I}_n^*$ then $p=2$ and $F_0$ depends on the component into
         which $\sigma_p$ specializes:
         $$\begin{array}{lll}
            \mbox{reduction type} & \mbox{specialization into} & F_0\\
             {\rm I}_0^*\mbox{ or } {\rm I}_1^* & & \mbox{1 point}\\
             {\rm I}_n^*, n\geq3, n \mbox{ odd} & \mbox{necessarily }\Theta_1 & \mbox{a curve}\\
             {\rm I}_n^*, n\geq2, n \mbox{ even} & \Theta_1 & \mbox{a curve}\\
                                                 & \Theta_2, \Theta_3 & \mbox{1 point,}
           \end{array}$$
         where the $\Theta_i$'s are those irreducible components of 
         multiplicity $1$ that do not intersect with $\sigma_0$. 
         Furthermore, $\Theta_2$ and $\Theta_3$ pass through the same 
         component of  multiplicity $2$.
  \end{enumerate}
\end{Proposition}

\proof
The generic fiber of the fibration is an ordinary elliptic curve, and $\sigma_p$
generates a subgroup scheme isomorphic to $\ZZ/p\ZZ$.
By \cite{Tate; Oort 1970}, this group scheme can either specialize to $\alpha_p$
or $\ZZ/p\ZZ$ in $X_0$.
Now, if $X_0$ is good and ordinary then the $p$-torsion subgroup scheme
$X_0[p]$ of $X_0$ is isomorphic to $(\ZZ/p\ZZ)\times\mu_p$,
which implies that $\sigma_p$ cannot meet $\sigma_0$ and $F_0=\emptyset$.
Similarly, if $X_0$ has multiplicative reduction then $\sigma_p$ has to
specialize into the component group and again $F_0=\emptyset$.
In case of good and supersingular reduction $X_0[p]$ is 
infinitesimal, which implies 
$(\sigma_0\cdot\sigma_p)_0\geq1$ and $F_0=X_0$.
In case of additive reduction and $(\sigma_0\cdot\sigma_p)_0\geq1$
then $\sigma_p$ induces an $\alpha_p$-action on $X_0$, thus
$F_0=X_0$.

We may thus assume that $X_0$ has additive reduction and
$(\sigma_0\cdot\sigma_p)_0=0$.
In particular, we obtain a non-trivial $\ZZ/p\ZZ$-action
on $X_0$.
Also, unless $X_0$ is of type ${\rm II}$, this fiber is a union
of $\PP^1$'s.
Moreover, $F_0$ is connected by \cite{Dolgachev; Keum 2001},
and is thus one point or a connected curve.
The next thing to note is that a $\ZZ/p\ZZ$-action on $\PP^1$ 
in characteristic $p$ has either precisely one fixed point
or the action is trivial.
Also, components of $X_0$ get mapped to components and
a point of $F_0$ where two components meet has
to be mapped to another such point under the $\ZZ/p\ZZ$-action.
From these facts one can easily work out $F_0$, which we
leave to the reader.
\qed\medskip

Let us recall from \cite{Dolgachev; Keum 2001} that the fixed
locus of a $\ZZ/p\ZZ$-action on a K3 surface is either
a finite set of at most two points or a connected curve.
Combining these results with 
Proposition \ref{fixed locus of translation} we obtain
our first structural result:

\begin{Theorem}
 \label{at most two fibers}
 Let $X\to\PP^1$ be an elliptic K3 surface with $p$-torsion section
 in positive characteristic $p$. 
 Then the fibration has at least one and at most two fibers that 
 are neither multiplicative nor ordinary.

 Moreover, if there are two such fibers then $p\leq5$, 
 these fibers have additive reduction,
 translation by $\sigma_p$ has precisely two fixed points and 
 the elliptic fibration
 arises as Frobenius pullback from an elliptic K3 surface.
\end{Theorem}

\proof
If the fibration has neither additive nor good supersingular fibers
then translation by $\sigma_p$ acts without fixed points by
Proposition \ref{fixed locus of translation}.
By \cite[Theorem 2.4]{Dolgachev; Keum 2001} this implies $p=2$
and that $Y=X/G$ is an Enriques surface, which
is absurd since genus one-fibrations on Enriques have
multiple fibers and thus can never be elliptic, i.e., with 
zero section.

The fixed locus of the $\sigma_p$-translation consists
either of at most two points or is a connected curve by
\cite{Dolgachev; Keum 2001}.
On the other hand, every fiber that is neither multiplicative
nor good ordinary has a non-trivial contribution to the
fixed locus by Proposition \ref{fixed locus of translation}.
This implies that there can be at most two fibers that
are good supersingular or additive.

Moreover, if there are two such fibers then the fixed locus
consists of two points.
By Proposition \ref{fixed locus of translation} these two 
fibers have additive reduction
and \cite[Theorem 2.4]{Dolgachev; Keum 2001}
implies $p\leq5$ and that $Y\to\PP^1$ is an elliptic K3 surface
(it cannot be Enriques by the reasons given above).
\qed\medskip

Next, the additive fibers tend to be potentially supersingular,
which is important for the computation of the formal Brauer
group in Section \ref{sec:brauer}.
The following extends results from
\cite{Liedtke; Schroeer 2008}.

\begin{Proposition}
  \label{potential reduction}
  Let $X\to B$ be an elliptic fibration with $p^n$-torsion sections
  and $p^n\geq3$.
  Then every additive fiber has potentially supersingular reduction.
\end{Proposition}

\proof
For $p\geq5$ this is \cite[Theorem 4.3]{Liedtke; Schroeer 2008}
and \cite[Remark 4.4]{Liedtke; Schroeer 2008}.

For $p^n=3$ and $p^n=4$ there is a universal elliptic curve over
$\Igord{p^n}$, which degenerates into multiplicative fibers at places
of potentially multiplicative reduction, see
\cite[Section 12]{Liedtke; Schroeer 2008} for $\Igord{3}$
and Proposition \ref{explicit equations} for $\Igord{4}$.
Since every elliptic fibration with $p^n$-torsion section pulls 
back from these, we conclude that the only additive fibers
can come from potentially supersingular places.
\qed

\begin{Proposition}
  \label{intersection numbers}
  Let $X\to\PP^1$ be an elliptic K3 surface with $p$-torsion section
  $\sigma_p$ in positive characteristic $p$.
  Then either
  \begin{enumerate}
    \item $\sigma_0\cdot\sigma_p=0$ and there are no fibers with
      good supersingular reduction, i.e., every potentially supersingular
      fiber has additive reduction, or
    \item $\sigma_0\cdot\sigma_p=1$, the characteristic is $p=2$,
       the fibration is semi-stable, and there is precisely
       one fiber with good supersingular reduction.
  \end{enumerate}
\end{Proposition}

\proof
Suppose that $\sigma_0\cdot\sigma_p\geq1$.
Then the fixed locus is a connected curve by 
\cite[Corollary 3.6]{Dolgachev; Keum 2001}
and from Proposition \ref{fixed locus of translation} and
Proposition \ref{at most two fibers} we infer that 
there is only one fiber whose reduction is neither good ordinary
nor multiplicative.
Moreover, the intersection of $\sigma_0$ and $\sigma_p$ takes place 
in this fiber.
In particular, there is at most one additive fiber and if
there is one, then $\sigma_p$ does not specialize into the
component group of that fiber.
We denote by ${\rm I}_{pn_v}$ with $v=1,...$ the multiplicative 
fibers and applying
\cite[Theorem 8.6]{Shioda 1990}, we get
$$ 
6\,\leq\,4+2(\sigma_0\cdot\sigma_p)\,=\,\sum_v n_v \frac{k_v(p-k_v)}{p},
$$
where the $k_v$'s are integers $1\leq k_v\leq p-1$ that encode which
component of ${\rm I}_{pn_v}$ is hit by $\sigma_p$.
Basic calculus tells us
$$
\frac{k_v(p-k_v)}{p}\,\leq\,\frac{p}{4}\,,
$$
where the inequality is strict if $p\neq2$.
On the other hand, we know
$$
24=c_2(X)=\sum_v pn_v + a
$$
where $a=0$ if and only if there are no additive fibers.
We conclude
$$
6\,\leq\,4+2(\sigma_0\cdot\sigma_p)\,\leq\,\sum_v n_v\frac{p}{4}\,\leq\,
\frac{c_2(X)}{4}\,=\,6,
$$
i.e, we have equality everywhere.
Thus, $\sigma_0\cdot\sigma_p=1$,
the characteristic equals $p=2$ (else the second inequality could not be an equality)
and there are no additive fibers (else the third inequality could
not be an equality).
\qed

\begin{Remark}
 The second alternative does exist and a complete classification is given 
 in Proposition \ref{p=2semistable}.
\end{Remark}

\section{Unirationality and the formal Brauer group}
\label{sec:brauer}

In this section we relate potentially supersingular fibres of an 
elliptic K3 surface with $p$-torsion sections to its
formal Brauer group.
For $p\geq3$ this implies that these surfaces
are either unirational or ordinary.
It also implies conjectures of Artin and Artin--Shioda
in this case.

We start by recalling the following fundamental result
of \cite{Artin; Mazur 1977}:
if $X$ is a smooth surface over $k=\overline{k}$ 
with smooth Picard scheme, e.g. a K3 surface,
then the functor on the category of finite local 
$k$-algebras $A$ with residue field $k$
$$\begin{array}{ccccc}
   \Brauer &:& A &\mapsto& 
   \ker\left(\Het{2}(X\times A,\GG_m)\to\Het{2}(X,\GG_m)\right)
\end{array}$$
is pro-represented by a smooth formal group of dimension $h^2(X,\OO_X)$,
the {\em formal Brauer group} $\Brauer(X)$ of $X$.

For a K3 surface, the height $h$ of the formal Brauer group 
is $\infty$ or an integer $1\leq h\leq10$ and all values are taken
\cite[Corollary 7.7]{Artin 1974}.
Moreover, $h$ determines the Newton polygon on second
crystalline cohomology
\cite[Section II.7.2]{Illusie 1979}.
In particular, the extreme cases are as follows:
\begin{itemize}
 \item[-] $h=1$ if and only if Newton- and
     Hodge- polygon coincide, i.e., the K3 surface is {\em ordinary}, and
 \item[-] $h=\infty$ if and only if the Newton polygon is a 
     straight line, i.e., the K3 surface is {\em supersingular}.
\end{itemize}

To be more precise about the notion of supersingularity, we recall

\begin{Definition}
 A K3 surface is called {\em supersingular in the sense of Artin}
 if its formal Brauer group has infinite height.
 A surface is called {\em supersingular in the sense of Shioda}
 if it satisfies $\rho=b_2$.
\end{Definition}

A K3 surface that is Shioda-supersingular is 
also Artin-supersingular \cite[Theorem 0.1]{Artin 1974}.
The Artin--Mazur conjecture states that also the converse holds
\cite[Remarque II.5.13]{Illusie 1979}.
Since this conjecture is known to be true for elliptic K3 surfaces
\cite[Theorem 1.7]{Artin 1974}, we do not
have to distinguish between these two notions 
of supersingularity.

Unirational K3 surfaces are Shioda-supersingular
\cite{Shioda 1974b}, as well as Artin-supersingular \cite{Artin 1974}.
For both notions, the converse is conjectured, see
\cite[Question II]{Shioda 1977} and \cite{Artin 1974}.
Thus, we summarize

\begin{Conjecture}
 \label{Artin Shioda Conjecture}
 For elliptic K3 surfaces,
 \begin{enumerate}
  \item (Shioda) Shioda-supersingularity implies unirationality,
  \item (Artin) Artin-supersingularity implies unirationality,
 \end{enumerate}
\end{Conjecture}

This conjecture is known to hold
for Fermat quartics \cite{Shioda 1974b},
Kummer surfaces in $p>2$ \cite{Shioda 1977b},
and thus for supersingular K3 surfaces
with Artin invariant $\sigma_0\leq2$ \cite{Ogus 1978}.
Also it holds in characteristic $2$ 
\cite{Rudakov Shafarevich 1979},
and for supersingular K3 surfaces with Artin invariant
$\sigma_0\leq3$ in characteristic $5$ \cite{Pho Shimada 2006}.

In characteristic $2$,
there is another conjecture by Artin
\cite[p.552]{Artin 1974},
which, if true, would imply the previous conjectures 
and gives a geometric explanation for them
-- note that this conjecture is supported by 
a dimension count \cite[p.552]{Artin 1974}:

\begin{Conjecture}[Artin]
 \label{Artin Conjecture2}
 In characteristic $2$, an elliptic fibration on a supersingular 
 K3 surface arises via Frobenius pullback from a rational elliptic surface.
\end{Conjecture}

Unfortunately, such a conjecture cannot be true in 
characteristic $p\geq3$.
Here is a counter-example:

\begin{Example}
Let $S_4$ be the Fermat quartic in $\PP^3$, which 
has been shown in \cite{Shioda 1974b}
to be supersingular in all characteristics $p$ for which
there exists a $\nu$ s.th. $p^\nu\equiv-1\mod4$, e.g.
in $p=3$.
This surface possesses a genus one fibration with six fibers
of type ${\rm I}_4$, see \cite[Section IV.2]{Barth; Hulek 1985}.
The associated Jacobian fibration $X\to\PP^1$ is 
a supersingular elliptic K3 surface, again with
six fibers of type ${\rm I}_4$. 
If it were the Frobenius pullback of some other elliptic
surface then the elliptic fibration of $X^{(1/p)}$ would have 
six fibers of type ${\rm I}_{n}$ such that
$pn=4$, giving $p=2$ as only possibility.
Thus, $X\to\PP^1$ is a supersingular K3 surface whose
elliptic fibration is not a Frobenius pullback from
another elliptic fibration.
\end{Example}

The following result links the height of the formal Brauer group to
the number of potentially supersingular fibers of the elliptic
fibration:

\begin{Theorem}
 \label{Brauer height}
 Let $X\to\PP^1$ be an elliptic K3 surface with $p$-torsion section
 in characteristic $p$, whose fibration does not have constant
 $j$-invariant. 
 Then the fibration has at least one and at most two fibres with 
 potentially supersingular reduction.
 Moreover, 
 \begin{enumerate}
  \item if there is one fiber with potentially supersingular reduction then
    the formal Brauer group has height $h\geq2$.
  \item if there are two fibers with potentially supersingular reduction then
    the formal Brauer group has height $h=1$.
 \end{enumerate}
\end{Theorem}

\prf
We know $p\leq7$ by Theorem \ref{possible cases}.
Since we assumed the fibration not to have constant $j$-invariant,
the map from the base to the $j$-line is dominant, 
whence surjective and there is at least one fiber with
potentially supersingular reduction.

Being a K3 surface, we may assume that the elliptic
fibration is given by a Weierstra\ss\ equation
\begin{equation}
\label{k3equation}
  y^2\,+\,a_1(t)xy\,+\,a_3(t)y\,=\,x^3\,+\,a_2(t)x^2\,+\,a_4(t)x\,+\,a_6(t)
\end{equation}
where the $a_i(t)$'s are polynomials of degree $\leq2i$, i.e.,
$a_i(t) \,=\,\sum_{j=0}^{2i} a_{ij}t^j$.

Assume $p=2$.
Then the formal Brauer group of $X$ has height $h=1$
if and only if $a_{11}\neq0$ by \cite[Theorem (2.12)]{Artin 1974}
(the extra assumptions of this theorem are not needed for this statement).
A fiber with potentially supersingular reduction is given by 
the vanishing of
$j(t)\,=\,a_1(t)^{12}/\Delta(t)$.
Since $\deg a_1(t)\leq2$, the fibration has at most two such fibers.
Moreover, if the fibration has two such fibres then $a_{11}\neq0$,
which implies $h=1$.
On the other hand, if there is only one such fiber then 
$a_{11}=0$, which implies $h\geq2$.

Now, assume that $p=3$.
A straight forward, but tedious calculation shows that 
$h=1$ is equivalent to $a_{11}^2+a_{22}\neq0$ in this case.
After a suitable change of coordinates, we may assume $a_1(t)=0$.
In this case, the Hasse invariant of the generic fiber is given
by the class of $-a_2(t)$ in $k(t)^\times/k(t)^{\times2}$.
Moreover, since the fibration has a $3$-torsion section,
the Hasse invariant is trivial, i.e., $-a_2(t)$ is a square.
On the other hand, fibers with potentially supersingular reduction
fulfill $0=c_4(t)=b_2(t)^2=a_2(t)^2$ in this case.
From $\deg a_2(t)\leq4$ we conclude that there are
at most two such fibers.
Moreover, the fibration has two such fibers if and only if 
$a_{22}\neq0$, i.e., if and only if $h=1$.

Next, assume $p=5$.
Then we may assume $a_1(t)=a_2(t)=a_3(t)=0$.
Computing the Hasse invariant, we see that then $2a_4(t)$ has to be
a fourth power in order for the fibration to possess
$5$-torsion sections.
The vanishing of $c_4(t)=2a_4(t)$ is necessary for a fiber
to have potentially supersingular reduction.
From $\deg a_4(t)\leq8$ and the fact that $2a_4(t)$ is a fourth
power we conclude
that there are at most two such fibers.
A tedious calculation shows that $h=1$ is equivalent to
$2a_{44}\neq0$ under our assumptions.
As in the previous cases, having two fibres with potentially
supersingular reduction is equivalent to $2a_{44}\neq0$, and thus 
equivalent to $h=1$.

We leave $p=7$ to the reader.
Alternatively, one can use 
Theorem \ref{p=7} below, by which there is only one
such surface.
It has one fiber with potentially supersingular reduction.
The elliptic fibration arises as Frobenius pullback from
a rational elliptic surface, i.e., this unique surface
is unirational, whence fulfills $h=\infty$.
\qed\medskip

\begin{Remark}
  In characteristic $2$,
  this connection between potentially supersingular
  fibers and the height of the formal Brauer group
  has already been observed by Artin \cite[p.552]{Artin 1974}.
\end{Remark}


Thus, in order to obtain
supersingular elliptic K3 surfaces with $p$-torsion sections
we have to look at fibrations with constant $j$-invariant, which can exist 
for $p=2$ only, or at fibrations that have precisely one potentially
supersingular fiber.

\begin{Proposition}
 \label{Frobenius from rational}
  Let $X\to\PP^1$ be an elliptic K3 surface with $p$-torsion section
  in characteristic $p\geq3$ that has precisely one potentially  
  supersingular fiber.
  Then the elliptic fibration arises as Frobenius pullback from
  a rational elliptic surface.
  In particular, $X$ is unirational and supersingular 
  ($h=\infty$).
\end{Proposition}

\proof
Let ${\rm I}_{pn_v}$, $v=1,...$ be the multiplicative fibers.
Since $p\geq3$, the fibration does not have constant $j$-invariant
and thus there exist places of potentially multiplicative reduction which
are multiplicative by Proposition \ref{potential reduction}.
Now, by Proposition \ref{intersection numbers} the potentially
supersingular fiber is additive, say with $m$ components
and Swan conductor $\delta$ and we obtain
\begin{equation}
\label{c2=24}
  24\,=\,c_2(X)\,=\,\sum_v pn_v\,+\,(2+\delta+(m-1))
\end{equation}
We also know that $X\to\PP^1$ arises as Frobenius pullback
from some elliptic fibration $Y\to\PP^1$, which 
has 
multiplicative fibers ${\rm I}_{n_v}$, $v=1,...$.
This fibration has one additive fiber also with
Swan conductor $\delta$ and with, say, $m'$ components.
Using (\ref{c2=24}) we obtain
\begin{equation}
 \label{c2estimate}
 c_2(Y)\,=\,\sum_v n_v + (2 + \delta + (m'-1)) \,\leq\,
 \frac{22-\delta}{p} + (2 + \delta + (m'-1))
\end{equation}
Since $p\neq2$, reduction of type ${\rm I}^*_n$ 
with $n\geq1$
is potentially multiplicative and thus cannot occur
as the additive fiber of $Y\to\PP^1$.
Inspecting the list of additive fibers we obtain $m'\leq9$.

On the other hand, $Y$ is either rational or K3, i.e,
$c_2(Y)=12$ or $c_2(Y)=24$.
If $p\geq5$ then $\delta=0$ and (\ref{c2estimate})
implies $c_2(Y)<24$, which implies that $Y$ is rational.
If $p=3$ then $c_2(Y)=24$ could only be achieved if
$\delta\geq20$.
However, this contradicts (\ref{c2=24}), since
$\sum_n pn_v\geq p=3$. 
Thus, $Y$ is a rational surface also for $p=3$.
\qed\medskip

\begin{Remark}
  We will see in Section \ref{sec:p=2} that 
  the statement is wrong for $p=2$.
\end{Remark}

We now come to one of the main results of this article,
which relates the geometry of the elliptic fibration
to supersingularity and unirationality.

\begin{Theorem}
 \label{supersingular characterisation}
 Let $X\to\PP^1$ be an elliptic K3 surface with $p$-torsion sections 
 in characteristic $p\geq3$.
 Let
 $\varphi:\PP^1\to\oIgord{p}$ be the compactified 
 classifying morphism.
 Then the following are equivalent:
 \begin{enumerate}
  \item $X$ arises as Frobenius pullback from a rational elliptic surface
  \item $X$ is a Zariski surface
  \item $X$ is unirational
  \item $X$ is supersingular
  \item the fibration has precisely one fiber with additive reduction
  \item $\varphi$ is totally ramified over the supersingular point of $\oIgord{p}$
 \end{enumerate}
 In particular, the conjectures of Artin--Shioda (Conjecture \ref{Artin Shioda Conjecture})
 and Artin (Conjecture \ref{Artin Conjecture2}) hold for this class of surfaces.
\end{Theorem}

\prf
Since $p\geq3$ the fibration does not have constant $j$-invariant 
by Theorem \ref{possible cases} and so $\varphi$ is surjective.
Also $p\leq7$ by loc. cit., which implies that $\oIgord{p}$ has precisely
one supersingular point.
By Proposition \ref{potential reduction}, all additive fibers are potentially
supersingular, which gives (5) $\Rightarrow$ (6).
By Proposition \ref{intersection numbers}, potentially supersingular fibers 
are additive and we get (6) $\Rightarrow$ (5).

The implications (1) $\Rightarrow$ (2) $\Rightarrow$ (3) $\Rightarrow$ (4) 
hold in general.
The implication (4) $\Rightarrow$ (5) follows from
Theorem \ref{Brauer height}. 
Finally,
(5) $\Rightarrow$ (1) follows from Proposition \ref{Frobenius from rational}.
\qed

\begin{Corollary}
 The Artin--Shioda conjecture holds for
 elliptic K3 surfaces with $p$-torsion sections.
\end{Corollary}

\prf
For $p\geq3$ this is Theorem \ref{supersingular characterisation}
and for $p=2$ it follows from \cite{Rudakov Shafarevich 1979}.
\qed

\begin{Remark}
 We will see in Section \ref{sec:p=2} that Theorem \ref{supersingular characterisation}
 also holds for elliptic K3 surfaces with $4$-torsion sections in 
 characteristic $2$.
\end{Remark}

Let us finally reformulate Theorem \ref{supersingular characterisation}
in terms of the "other" surfaces:

\begin{Theorem}
 \label{ordinary characterisation}
 Let $X\to\PP^1$ be an elliptic K3 surface with $p$-torsion sections 
 in characteristic $p\geq3$.
 Then the following are equivalent:
 \begin{enumerate}
  \item $X$ is ordinary
  \item $X$ is not unirational
  \item $X$ arises as Frobenius pullback from a K3 surface
  \item the fibration has precisely two fibers with additive reduction
 \end{enumerate}
 Moreover, such surfaces can exist in characteristic $p\leq5$ only.
\end{Theorem}

\proof
By Theorem \ref{at most two fibers} case (4) can happen in
characteristic $p\leq5$ only.

The implications (1) $\Rightarrow$ (2) $\Rightarrow$ (3) hold in 
general.
The implication (3) $\Rightarrow$ (4) follows from 
Theorem \ref{supersingular characterisation} and
Theorem \ref{at most two fibers}.
Finally, the implication (4) $\Rightarrow$ (1) follows
from Theorem \ref{Brauer height}.
\qed

\begin{Remark}
 As we shall see in Section \ref{sec:p=2}, there do exist elliptic
 K3 surfaces with $2$-torsion section that are neither
 unirational nor ordinary.
\end{Remark}

\section{The explicit classification}
\label{sec:classification}

Having established the general picture in the 
previous sections,
we now give a detailed classification of elliptic
K3 surfaces with $p$-torsion section in characteristic
$p\geq3$.
This is achieved by studying the classifying morphism
to the Igusa curve and the N\'eron model of the universal
family over $\Igord{p}$.
We pay special attention to the arising supersingular
surfaces.

Let us recall the following from \cite{Shioda 1990}:
the N\'eron--Severi group $\NS(X)$ of an elliptic
surface together with its intersection pairing
is made up of two natural subgroups:
the trivial lattice $T$, which is associated to 
the singular fibers, and the 
{\em Mordell--Weil group}
$\MW(X)$, which arises from sections of the elliptic
fibration and the N\'eron--Tate height pairing.
Inside this group sits the {\em narrow Mordell--Weil group}
$\nMW(X)$ consisting of those sections 
that lie fiberwise on the same component as the
zero-section.
For rational elliptic surfaces these groups
have been worked out explicitly
in \cite{Oguiso; Shioda 1991}.

For the singular fibers we use Kodaira's notation.
For example, ${\rm I}_n$ denotes the multiplicative
reduction where a singular fiber consists of 
$n$ smooth rational curves forming a cycle.
In case of additive reduction and in characteristic
$p\leq3$ there is a further invariant, namely the
{\em Swan conductor} $\delta$ of a singular fiber,
which we add as index.
Thus, ${\rm I}_{n,\delta}^*$ stands for additive
reduction of type ${\rm I}_n^*$ with Swan conductor
$\delta$.
We refer to \cite[Chapter IV]{Silverman 1994} for 
definitions and details.

Let us also recall that the discriminant
of $\NS(X)$ for a supersingular K3 surface $X$
is of the form $p^{2\sigma_0}$ for some integer
$1\leq\sigma_0\leq 10$, called the {\em Artin invariant}
\cite{Artin 1974}.
All values actually occur and surfaces with Artin
invariant $\leq\sigma_0$ form a $(\sigma_0-1)$-dimensional
subspace inside the moduli space of all
supersingular surfaces \cite{Ogus 1978}.
Finally, there is only surface with $\sigma_0=1$
\cite{Ogus 1978} and in $p\geq3$ surfaces with
$\sigma_0\leq2$ are Kummer surfaces by
\cite{Ogus 1978} and \cite{Shioda 1979}.

\subsection*{Characteristic $\mathbf{7}$}

\begin{Theorem}
  \label{p=7}
  There exists only one elliptic K3 surface 
  $X\to\PP^1$ with $7$-torsion section in characteristic $7$
  up to isomorphism.
  It has the following invariants:
  $$\begin{array}{lccc}
     \mbox{singular fibers} & \sigma_0 & \nMW(X) & \MW(X) \\
     \hline
      {\rm III}, 3\times{\rm I}_7 & 1 & A_1(7) & A_1^*(7)\oplus(\ZZ/7\ZZ) 
    \end{array}
  $$
  The Weierstra\ss\ equation is given by the following:
  $$ 
     y^2\,=\,x^3\,+\,tx\,+\,t^{12}.
  $$
  In particular, it is the unique supersingular K3 surface with
  Artin invariant $\sigma_0=1$.
\end{Theorem}

\prf
As in Section \ref{sec:generalities},
we denote by $\varphi$ the classifying morphism to $\oIgord{7}$
and by ${\cal E}\to\Igord{7}$ the universal curve.
An analysis of the multiplicative fibers 
as in the proof of Theorem \ref{possible cases} shows
that $\deg\varphi\geq2$ is impossible.
Hence $\varphi$ is an isomorphism, proving uniqueness.
Since ${\cal E}^{(7)}$ corresponds in fact a K3 surface,
we get existence.
The singular fibres are listed in 
Proposition \ref{explicit equations}.

Denote by $Y\to\PP^1$ the elliptic fibration corresponding to 
${\cal E}$.
Then $Y$ is rational, which implies that
$X$ is a Zariski surface and thus unirational.
The singular fibres are given in Proposition \ref{explicit equations}
and thus the root lattice of $Y$ is $E_7$.
From the tables in \cite{Oguiso; Shioda 1991} we see that 
the (narrow) Mordell-Weil lattice is 
$\MW(Y)\iso A_1^*$ and
$\nMW(Y)\iso A_1$, respectively.

Now, Frobenius induces an incluses of lattices
\begin{equation}
 \label{lattice inclusion}
 {\rm MW(Y)}_{\rm free}(p) \,\subseteq {\rm MW(X)}_{\rm free},
\end{equation}
which is of some finite index $\mu$.
Taking determinants, we obtain
$$
   \mu^2\,=\,\frac{\det\MW(Y)_{\rm free}(p)}{\det\MW(X)_{\rm free}}.
$$
After plugging in Lemma \ref{determinant formula} below, we obtain
$$
\mu^2\,=\,\frac{\det A_1^*(7)}{\det\NS(X)\,|\MW(X)_{\rm tor}|^2}\det(U\oplus A_6^{\oplus3}\oplus A_1)
\,=\,\frac{\frac{1}{2}\cdot7}{7^{2\sigma_0(X)}\,\cdot\,7^2}\cdot7^3\cdot2,
$$
which yields $\mu=1$. 
Thus, $\sigma_0=1$ and $\MW(X)\iso A_1^*(7)\oplus(\ZZ/7\ZZ)$.
\qed

\begin{Remark}
 Existence and uniqueness of this surface have already been shown in
 \cite[Examples 2.4]{Schweizer 2005}.
\end{Remark}

\begin{Lemma}[{\cite[Theorem 8.7]{Shioda 1990}}]
  \label{determinant formula}
  Let $X$ be an elliptic surface whose $j$-invariant is not constant. Then 
  $$
   \det\NS(X)\,=\,\frac{\det\MW(X)_{\rm free}\cdot\det T}{|\MW(X)_{\rm tor}|^2}\,,
  $$
  where $T$ denotes the trivial lattice.
  \qed
\end{Lemma}

\subsection*{Characteristic $\mathbf{5}$}

\begin{Theorem}
 \label{p=5}
 In characteristic $5$, the classifying morphism $\varphi$
 of an elliptic K3 surface with
 $5$-torsion section is finite of degree $2$.
 Conversely, if $\varphi:\PP^1\to\oIgord{5}$ 
 is a morphism of degree $2$ 
 then the associated elliptic fibration with 
 $5$-torsion section is a K3 surface.

 More precisely, the surfaces have the following invariants:
 $$\begin{array}{llcccc}
   & \mbox{singular fibers} & \dim & \sigma_0 & \nMW(X) & \MW(X) \\
 \hline
   &  2\times{\rm II}, 4\times{\rm I}_5 & 2 &  \\
   &  2\times{\rm II}, {\rm I}_{10}, 2\times{\rm I}_5 & 1 & \\
   &  2\times{\rm II}, 2\times{\rm I}_{10} & 0 & \\
   &  {\rm IV}, 4\times{\rm I}_5 & 1 & 2 & A_2(5) & A_2^*(5)\oplus\ZZ/5\ZZ \\
   &  {\rm IV}, {\rm I}_{10}, 2\times{\rm I}_5 & 0 & 1 & \langle 30 \rangle &
       \langle\frac{5}{6}\rangle\oplus\ZZ/5\ZZ 
 \end{array}$$
 Here, $\dim$ denotes the dimension of the family.
 For the supersingular surfaces, this list also gives 
 Artin invariants $\sigma_0$ and their (narrow) Mordell--Weil 
 lattices.
\end{Theorem}

\begin{Remark}
 \label{rem:supersing5}
 The surfaces with two ${\rm II}$-fibers arise as Frobenius pullbacks
 from Shioda's sandwich surfaces \cite{Shioda 2006}. 
 From this fact one obtains another proof of their 
 non-supersingularity.
\end{Remark}

\proof 
The proof is analogous to the proof of Theorem \ref{p=7}.
We leave it to the reader to show that the classifying
morphism $\varphi$ is of degree $2$.
Then we obtain the complete classification of these
surfaces in terms of the branch points of 
the classifying morphism:
To do so, let ${\cal E}\to\oIgord{5}$ be the universal elliptic
curve over the Igusa curve.
By Proposition \ref{explicit equations} its
Weierstra\ss\ equation is given by
$$
  y^2\,=\,x^3\,+\,3t^4x\,+\,t^5\,,
$$
which has a singular fiber of type ${\rm II}^*$ over $t=0$
and fibers of type ${\rm I}_1$ over $t=\pm 1$.
Note that this surface is a rational extremal elliptic surface.

We write the classifying morphism 
$\varphi=\varphi_{\alpha \beta}:\PP^1\to\oIgord{5}$ as 
$$
  t\,=\, \frac{\alpha s^2\,+\,\beta }{s^2\,+\,1}\,
$$
whose branch points are $t=\alpha$ and $t=\beta$, where
$t$ (resp. $s$) is a local parameter of $\oIgord{5}$
(resp. $\PP^1$).
Then our surfaces arise as pull-backs along Frobenius $F$ and $\varphi_{\alpha \beta}$:
$$
\begin{CD}
   X = Y^{(p)} @>>> Y @>>> {\cal E} \\
   @VVV @VVV  @VVV \\
   \PP^1 @>{F}>> \PP^1 @>\varphi_{\alpha \beta}>> 
   \oIgord{5}
\end{CD}
$$
The elliptic surface $Y$ is given by the Weierstra\ss\
equation
$$
  y^2\,=\,x^3\,+\,3(\alpha s^2+\beta)^4x\,+\,(\alpha s^2+\beta)^5\,(s^2+1),
$$
and depending on $\alpha$ and $\beta$ we obtain the following list
$$\begin{array}{lllc}
   \{\alpha, \beta\}\cap \{0, \pm 1\} 
   & \mbox{singular fibers of $X$} 
   & \mbox{singular fibers of $Y$} 
   & Y \\
\hline
   \emptyset
   & 2\times {\rm II}, 4\times {\rm I}_5 
   & 2\times {\rm II}^*, 4\times {\rm I}_1
   & \mbox{K3} \\
   \{1\},  \{-1\} 
   & 2\times {\rm II}, {\rm I}_{10}, 2\times {\rm I}_5 
   & 2\times {\rm II}^*, {\rm I}_2, 2\times {\rm I}_1
   & \mbox{K3} \\
   \{1, -1\} 
   & 2\times {\rm II}, 2\times {\rm I}_{10} 
   & 2\times {\rm II}^*, 2\times {\rm I}_2
   & \mbox{K3} \\
   \{0\}
   & {\rm IV}, 4\times {\rm I}_5 
   & {\rm IV}^*, 4\times {\rm I}_1
   & \mbox{rational} \\
   \{0, 1\}, \{0,-1\} 
   & {\rm IV}, {\rm I}_{10}, 2\times {\rm I}_{5} 
   & {\rm IV}^*, {\rm I}_2, 2\times {\rm I}_1 
   & \mbox{rational} \\
\end{array}$$
giving the explicit classification of our surfaces.

By Theorem \ref{supersingular characterisation} the supersingular surfaces
are precisely those that arise as Frobenius pullbacks from rational
elliptic surfaces.
It remains to determine the Mordell--Weil groups and Artin invariants.

For the $({\rm IV},{\rm I}_{10}, 2\times{\rm I}_5)$-surface
this can be done as in the proof of Theorem \ref{p=7} and we leave
it to the reader.

Let $X\to\PP^1$ be a $({\rm IV}, 4\times{\rm I}_5)$-surface.
Using \cite{Oguiso; Shioda 1991}, we see that it
arises via Frobenius pullback from rational
elliptic surface  $Y\to\PP^1$ with $\nMW(Y)\iso A_2$.
From (\ref{lattice inclusion}) we get an inclusion of
Mordell--Weil lattices and once we have shown equality
our assertion follows.
Now, $\nMW(Y)$ is generated by two sections $P_1,P_2$
with $\langle P_i,P_i\rangle=2$, which implies that both neither
meet the zero-section nor specialize into the component groups of the
singular fibers.
By Lemma \ref{verschiebung} below, these two sections cannot lie
in the image of $V:\MW(X)\to\MW(Y)$.

Now, denote by $K$ the function field of $\PP^1$
and let $E$ and $E^{(p)}$ be the generic fibers of $Y$
and $X$ over $\Spec K$.
Multiplication by $p$ induces an exact sequence
\begin{equation}
\label{VFpsequence}
\begin{array}{cccccccc}
0&\to&\ker(V)&\to&E^{(p)}(K)/F(E(K))&\stackrel{\overline{V}}{\to}&
E(K)/pE(K)\\
&&&\to&E(K)/V(E^{(p)}(K))&\to&0&,
\end{array}
\end{equation}
where $V$ denotes Verschiebung.
Knowing that $P_1$ and $P_2$ do not lie in the image
of $V$, this implies $\overline{V}=0$ in the sequence above
and we obtain the desired equality of Mordell--Weil lattices.
\qed

\begin{Lemma}
  \label{verschiebung}
  Let $R$ be complete DVR with field of fractions $K$ of 
  characteristic $p\geq5$ and perfect residue field $k$.
  Let $E$ be an elliptic curve over $K$ and assume that 
  $E^{(p)}$ has a $K$-rational $p$-division point. 
  Assume moreover, that $E$ has additive reduction that
  is not of type ${\rm II}^*$ if $p=5$.
  If $P\in E^{(p)}(K)$ then $V(P)$, where $V$ denotes Verschiebung, 
  specializes into the component group or to zero in 
  the N\'eron model ${\cal E}$ of $E$.
\end{Lemma}

\prf
Let $\pi\in R$ be a uniformizer.
Set $L:=K(\pi^{1/12})$ and denote by $S$ the integral closure of $R$
in $L$.
Then $L/K$ is totally ramified, $\varpi:=\pi^{1/12}$ is a uniformizer
on $S$.
Denote by $\nu_\pi$ and $\nu_\varpi$ normalized valuations, i.e.,
$\nu_\pi(\pi)=\nu_\varpi(\varpi)=1$ and
$\nu_\varpi(x)=12\nu_\pi(x)$ for all $x\in R$.

Since $p\geq5$, the curve $E$ acquires semi-stable reduction 
over $L$, which is good and supersingular 
\cite[Theorem 4.3]{Liedtke; Schroeer 2008}.
Let us denote by $\cal E$ minimal Weierstra\ss\ equations
and assume that the singularity (in case of bad reduction) 
lies in $(0,0)$.

For a section $P=(x_{0,K},y_{0,K})$ we set 
$t_{0,K}:=y_{0,K}/x_{0,K}$ and note that
$\nu_\pi(t_{0,K})<0$ if and only if $P$ specializes to zero
in the N\'eron model, as well as $\nu_\pi(t_{0,K})>0$ 
if and only if $P$
specializes non-trivially into the component group.
Now, we run Tate's algorithm and
suppose we have to reduce $r_1$-times to get
from ${\cal E}_K^{(p)}\times_K L$ to ${\cal E}_L^{(p)}$.
By our assumptions on $p$ and $L/K$ 
we have $r_1=\nu_\pi(\Delta_{\rm min})$, 
where $\Delta_{\rm min}$ denotes the minimal
discriminant of ${\cal E}_K^{(p)}$.
Then $P$, considered as a section of ${\cal E}_L^{(p)}$, fulfills
$\nu_{\varpi}(t_{0,L})=12\nu_{\pi}(t_{0,K})-r_1$.

Next, $V$ induces a map ${\cal E}_L^{(p)}\to {\cal E}_L$.
Both elliptic curves have good supersingular reduction and on the
level of tangent spaces, this map is multiplication by the
Hasse invariant \cite[Chapter 12.4]{Katz; Mazur 1985}.
Then, for appropriate local parameters $e$, $e^{(p)}$ around 
zero, $V$ is given by $e^{(p)}\mapsto H\cdot e+...$ for 
some lift of the Hasse invariant to $S$.
As this lift we may choose the ''naive`` Hasse invariant 
in the sense of raising a homogeneous Weierstra\ss\ equation to the 
$(p-1)$.st power and taking the coefficient of $(xyz)^{p-1}$.
If we set $h:=\nu_\varpi(H)$ then, $h>0$ since we have supersingular
reduction and $h$ is divisible by $p-1$ as there is
an $L$-rational $p$-division point on $E_L$.
Thus, if $V(P)=(x_{0,L}',y_{0,L}')$ in ${\cal E}_L$, 
we set $t_{0,L}'=y_{0,L}'/x_{0,L}'$ and get
$\nu_\varpi(t_{0,L}')= (12\nu_\pi(t_{0,K})-r_1)(h+1)$.

Suppose we have to reduce $r_2$-times in the Tate algorithm 
to get from ${\cal E}_K\times L$ to ${\cal E}_L$.
Then we finally obtain
$$
 \nu_\pi(t_{0.K}')\,=\,
 \nu_\pi(t_{0,K})+\frac{r_2-r_1}{12}+\left(\nu_\pi(t_{0,K})-\frac{r_1}{12}\right)\cdot h
$$
Let us first assume that $P\in E_K(K)$ does not 
specialize into the component group, which means $\nu_\pi(t_{0,K})\leq0$.
Recall that $h>0$ and that $p-1$ divides $h$.
Moreover, from the tables of minimal discriminants we get $r_2-r_1\leq8$
(note that reduction of type ${\rm I}^*_n$, $n\geq2$ is impossible
by \cite[Corollary 4.5]{Liedtke; Schroeer 2008}).
Thus, if $p\geq7$ or if $r_2-r_1<8$ we get $\nu_\pi(t_0'')<0$, i.e.,
$V(P)$ specializes to zero in the N\'eron model of $E_K$.
The only case where this may fail is $p=5$ and $r_2-r_1=8$, i.e.,
${\cal E}_K^{(p)}$ has reduction of type ${\rm II}$ ($r_1=2$)
and 
${\cal E}_K$ has reduction of type ${\rm II}^*$ ($r_2=10$).

Finally, assume that $P$ specializes into the component group of ${\cal E}_K^{(p)}$.
Then there exists an integer $m$, prime to $p$, such that $mP$ does {\em not}
specialize into the component group any more.
By the previous discussion $V(mP)$ specializes to zero in the N\'eron
model of ${\cal E}_K$.
Now, as a group scheme, the special fiber of ${\cal E}_K$ is 
$\GG_a\times\Phi$,
where $\Phi$ is the component group of ${\cal E}_K$.
Since $\GG_a$ does not have $m$-torsion, it follows that 
$V(P)$ specializes to zero or into the component group of
${\cal E}_K$.
\qed\medskip

The following result makes sure that we find in fact complete
families of supersingular K3 surfaces.

\begin{Proposition}
 \label{complete families}
 Let $X$ be an elliptic K3 surface with $p^n$-torsion section
 in characteristic $p$.
 Assume that $X$ is supersingular with Artin-invariant
 $\sigma_0$.
 Then, every (Shioda-)supersingular K3 surface with 
 Artin invariant $\sigma_0$ in characteristic $p$ 
 possesses an elliptic fibration with $p^n$-torsion section. 
\end{Proposition}

\prf
To give a (quasi-)elliptic fibration on $X$ is equivalent
to giving an isometric embedding of a hyperbolic
lattice $U$ of rank $2$ into $\NS(X)$.

Then, the trivial lattice $T$ is the sub-lattice
of $\NS(X)$ generated by $U$ and all $x\in U^\perp$
with $x^2=-2$, see \cite{Shioda 1990}.
By \cite[Theorem 1.3]{Shioda 1990} the torsion sections
of the fibration correspond to the torsion
of $\NS(X)/T$.

The N\'eron--Severi group of a (Shioda-)supersingular 
K3 surface is uniquely determined by $p$ and $\sigma_0$ by
\cite[Theorem 2']{Rudakov Shafarevich 1979}.
Thus, by the previous discussion, if one of these surfaces 
possesses a (quasi-)elliptic fibration with $p^n$-torsion section
then so do all of them.

However, we have to rule out the possibility that the
isometric embedding of $U$ into $\NS(X)$ 
corresponding to the elliptic fibration on $X$ gives rise
to a quasi-elliptic fibration on another K3 surface $Y$ 
with the same $p$ and $\sigma_0$:
if $p\geq5$ or if ${\rm rank}(T)<22$ then the fibration 
on $Y$ is automatically elliptic and the quasi-elliptic case
cannot occur at all.
And finally,
if $p\leq3$ and ${\rm rank}(T)=22$ then the elliptic fibration
on $X$ is extremal and these K3 surfaces have been
explicitly classified in \cite{Ito 2002}.
It turns out that these surfaces have Artin invariant
$\sigma_0=1$, i.e., $X$ is isomorphic to $Y$. 
\qed\medskip

Together with Theorem \ref{p=5} we immediately conclude

\begin{Corollary}
 Every (Shioda-)supersingular K3 surface with $\sigma_0\leq2$ 
 in characteristic $5$ possesses an elliptic
 fibration with $5$-torsion section.\qed
\end{Corollary}

\subsection*{Characteristic $\mathbf{3}$}

We denote by $O\in\oIgord{3}$ the unique supersingular point.

\begin{Theorem}
 \label{p=3}
 In characteristic $3$, the classifying morphism $\varphi$
 for an elliptic K3 surface with
 $3$-torsion section is finite of degree fulfills 
 $2\leq\deg\varphi\leq6$.
 More precisely,
 \begin{enumerate}
  \item $\deg\varphi=2$ and $\varphi^{-1}(O)$ consists of
     two points.
  \item $\deg\varphi=3$, $\varphi$ is separable and
     $\varphi^{-1}(O)$ consists of two points.
  \item $\deg\varphi=4$ and $\varphi^{-1}(O)$ consists of
     one or two points.
  \item $\deg\varphi=5$ and $\varphi^{-1}(O)$ consists of
     one point or two points with ramification index 
     $e=2$ and $e=3$.
  \item $\deg\varphi=6$ and $\varphi^{-1}(O)$ consists of
     one point or two points with ramification index 
     $e=3$.
 \end{enumerate}
 Conversely, if $\varphi$ is as above 
 then the associated elliptic fibration with 
 $3$-torsion section is a K3 surface.
\end{Theorem}

\proof
The proof is analogous to the proof of Theorem \ref{p=7}
(but lengthier and with more subcases)
and we leave it to the reader.
\qed\medskip

From this description it is easy to obtain a complete list of these
surfaces as before.
However, since this list is rather long, 
we have decided not to include it here.
Instead, we only determine the supersingular K3 surfaces 
with $3$-torsion sections.
By Theorem \ref{supersingular characterisation},
these are precisely the surfaces, where the
classifying morphism is totally ramified over
$O\in\oIgord{3}$.
As before, $\varphi$ denotes the classifying
morphism.

\begin{Theorem}
 \label{char. 3 table}
 Every (Shioda-)supersingular K3 surface with Artin invariant
 $\sigma_0\leq6$ in characteristic $3$ possesses an elliptic fibration with
 $3$-torsion section. 

 The complete list of these surfaces
 is given by the following table:
\medskip

 $\deg \varphi = 6$ (separable)
$$\begin{array}{lcccc}
         \mbox{singular fibers}
         & \dim         
         & \sigma_0
         & \nMW(X) 
         & \MW(X) 
  \\
   \hline
         {\rm II}_4, 6\times {\rm I}_3
         & 5 
         & 6 
         & E_8(3)
         & E_8(3)\oplus {\ZZ/3\ZZ} 
  \\
         {\rm II}_4, {\rm I}_6, {\rm I}_3 \times 4
         & 4 
         & 5 
         & E_7(3)
         & E_7^*(3)\oplus {\ZZ/3\ZZ} 
  \\
         {\rm II}_4 , {\rm I}_9, {\rm I}_3 \times 3
         & 3 
         & 4 
         & E_6(3)
         & E_6^*(3)\oplus {\ZZ/3\ZZ}
  \\
         {\rm II}_4, {\rm I}_6 \times 2, {\rm I}_3 \times 2
         & 3 
         & 4 
         & D_6(3)
         & D_6^*(3)\oplus {\ZZ/3\ZZ}
  \\
         {\rm II}_4 , {\rm I}_{12}, {\rm I}_3 \times 2 
         & 2
         & 3
         & D_5(3)
         & D_5^*(3)\oplus {\ZZ/3\ZZ}
  \\
         {\rm II}_4, {\rm I}_6 \times 3
         & 2 
         & 3 
         & D_4(3)\oplus A_1(3)
         & D_4^*(3)\oplus A_1^*(3)\oplus {\ZZ/3\ZZ}
  \\
         {\rm II}_4, {\rm I}_9, {\rm I}_6, {\rm I}_3
         & 2 
         & 3 
         & A_5(3)
         & A_5^*(3)\oplus {\ZZ/3\ZZ}
  \\
         {\rm II}_4, {\rm I}_{15}, {\rm I}_3
         & 1
         & 2
         & A_4(3)
         & A_4^*(3)\oplus {\ZZ/3\ZZ}
  \\
         {\rm II}_4, {\rm I}_{12}, {\rm I}_6
         & 1
         & 2
         & A_3(3)\oplus A_1(3)
         & A_3^*(3)\oplus A_1^*(3)\oplus {\ZZ/3\ZZ}
  \\
         {\rm IV}_2,  6\times {\rm I}_3
         & 4
         & 5 
         & E_6(3)
         & E_6^*(3)\oplus {\ZZ/3\ZZ}
  \\
         {\rm IV}_2, {\rm I}_6, {\rm I}_3 \times 4
         & 3 
         & 4 
         & A_5(3)
         & A_5^*(3)\oplus {\ZZ/3\ZZ}
  \\
         {\rm IV}_2 , {\rm I}_9, {\rm I}_3 \times 3
         & 2 
         & 3 
         & A_2(3)^{\oplus 2}
         & {A_2^*}(3)^{\oplus 2}\oplus {\ZZ/3\ZZ}
  \\
         {\rm IV}_2, {\rm I}_6 \times 2, {\rm I}_3 \times 2
         & 2 
         & 3 
         & L_4(3)
         & L_4^*(3)\oplus {\ZZ/3\ZZ}
  \\
         {\rm IV}_2 , {\rm I}_{12}, {\rm I}_3 \times 2
         & 1 
         & 2 
         & L_3(3)
         & L_3^*(3)\oplus {\ZZ/3\ZZ}
  \\
         {\rm IV}_2, {\rm I}_6 \times 3
         & 1 
         & 2 
         & A_1(3)\oplus L_2(3)
         & A_1^*(3)\oplus L_2^*(3)\oplus {\ZZ/3\ZZ}
  \\
         {\rm I}_{0,0}^*, {\rm I}_3 \times 6
         & 3
         & 4 
         & D_4(3)
         & D_4^*(3)\oplus {\ZZ/3\ZZ}
  \\
         {\rm I}_{0,0}^*, {\rm I}_6, 4 \times {\rm I}_3 
         & 2
         & 3 
         & A_1(3)^{\oplus 3}
         & {A_1^*(3)}^{\oplus 3}\oplus {\ZZ/3\ZZ} 
  \\
         {\rm I}_{0,0}^*, {\rm I}_6 \times 2, {\rm I}_3 \times 2
         & 1 
         & 2 
         & A_1(3)^{\oplus 2}
         & {A_1^*}(3)^{\oplus 2}\oplus \ZZ /6\ZZ
  \\
         {\rm I}_{0,0}^*, {\rm I}_9, {\rm I}_3 \times 3
         & 1 
         & 2 
         & L_2(3) 
         & L_2^*(3)\oplus {\ZZ/3\ZZ} 
  \\
         {\rm I}_{0,0}^*, {\rm I}_6 \times 3
         & 0 
         & 1 
         & A_1(3)
         & A_1^*(3)\oplus \ZZ / 6\ZZ\oplus {\ZZ/2\ZZ}
  \\
         {\rm I}_{0,0}^*, {\rm I}_{12}, {\rm I}_3 \times 2 
         & 0 
         & 1 
         & \langle 12\rangle 
         & \langle \frac{3}{4} \rangle \oplus \ZZ/ 6\ZZ
\end{array}$$

 $\deg \varphi = 6$ (inseparable)
$$\begin{array}{lcccc}
         \mbox{singular fibers}
         & \dim         
         & \sigma_0
         & \nMW(X) 
         & \MW(X) 
  \\
\hline   
         {\rm IV}_2, 2\times {\rm I}_9 
         & 1 
         & 2 
         & A_2(3)
         & A_2^*(3)\oplus {\ZZ/3\ZZ}
  \\
         {\rm IV}_2, {\rm I}_{18}
         & 0 
         & 1 
         & \langle 18\rangle
         & \langle \frac{1}{2}\rangle\oplus {\ZZ/3\ZZ}
  \\
\end{array}$$

 $\deg \varphi = 5$
$$\begin{array}{llcccc}
         \mbox{singular fibers}
         & \dim         
         & \sigma_0
         & \nMW(X) 
         & \MW(X) 
  \\
\hline
         {\rm IV}_5, 5\times {\rm I}_3
         & 4
         & 5 
         & E_8(3)
         & 3.(E_8(3))\oplus {\ZZ/3\ZZ}
  \\
         {\rm IV}_5, {\rm I}_6, 3\times {\rm I}_3
         & 3 
         & 4 
         & E_7(3)
         & 3.(E_7^*(3))\oplus {\ZZ/3\ZZ}
  \\
         {\rm IV}_5, 2\times {\rm I}_6 , {\rm I}_3
         & 2 
         & 3 
         & D_6(3)
         & 3.(D_6^*(3))\oplus {\ZZ/3\ZZ}
  \\
         {\rm IV}_5, {\rm I}_9, 2 \times {\rm I}_3
         & 2 
         & 3 
         & E_6(3)
         & 3.(E_6^*(3))\oplus {\ZZ/3\ZZ}
  \\
         {\rm IV}_5, {\rm I}_9, {\rm I}_6
         & 1 
         & 2 
         & A_5(3)
         & 3.(A_5^*(3))\oplus {\ZZ/3\ZZ}
  \\
         {\rm IV}_5, {\rm I}_{12}, {\rm I}_3
         & 1 
         & 2 
         & D_5(3)
         & 3.(D_5^*(3))\oplus {\ZZ/3\ZZ}
  \\
         {\rm IV}_5, {\rm I}_{15}
         & 0 
         & 1 
         & A_4(3)
         & 3.(A_4^*(3))\oplus {\ZZ/3\ZZ}
\end{array}$$

 $\deg \varphi = 4$
$$\begin{array}{llcccc}
         \mbox{singular fibers}
         & \dim         
         & \sigma_0
         & \nMW(X) 
         & \MW(X) 
  \\
 \hline
         {\rm IV}_4^*, 4\times {\rm I}_3 
         & 3
         & 4  
         & E_6(3)
         & E_6^*(3)\oplus {\ZZ/3\ZZ}
  \\
         {\rm IV}_4^*, {\rm I}_6, 2\times {\rm I}_3
         & 2  
         & 3 
         & A_5(3)
         & A_5^*(3)\oplus {\ZZ/3\ZZ}
  \\
         {\rm IV}_4^*, {\rm I}_6, {\rm I}_6 
         & 1 
         & 2 
         & L_4(3)
         & L_4^*(3)\oplus {\ZZ/3\ZZ}
  \\
         {\rm IV}_4^*, {\rm I}_{9}, {\rm I}_3 
         & 1 
         & 2 
         & A_2(3)^{\oplus 2}
         & {A_2^*}(3)^{\oplus 2}\oplus {\ZZ/3\ZZ}
  \\
         {\rm IV}_4^*, {\rm I}_{12}
         & 0  
         & 1  
         & L_3(3)
         & L_3^*(3)\oplus {\ZZ/3\ZZ}
\end{array}$$
Here, $L_2$, $L_3$, and $L_4$ are lattices of rank $2$, $3$, $4$, all of
determinant $12$, whose matrices are given by 
$$L_2\,=\,\begin{pmatrix}
4 & -2 \\ -2 & 4
\end{pmatrix},\,
L_3\,=\,\begin{pmatrix}
2 & 0 & -1 \\
0 & 2 & -1 \\
-1 & -1 & 4 
\end{pmatrix},\, 
L_4\,=\,\begin{pmatrix}
4 & -1 & 0 & 1 \\
-1 & 2 & -1 & 0 \\
0 & -1 & 2 & -1 \\
1 & 0 & -1 & 2
\end{pmatrix}.$$
Also, the notation $3.L$ for a lattice $L$ stands for a lattice that
has $L$ as a sublattice of index $3$.
%
%
\end{Theorem}


\prf
By Theorem \ref{supersingular characterisation} the classifying
morphism $\varphi$ is totally ramified over the supersingular
point $O\in\oIgord{3}$.
This gives $4\leq\deg\varphi\leq6$ by Theorem \ref{p=3}.
We proceed as in the proof of Theorem \ref{p=5}
in order to obtain explicit equations:
let $f_3(s)$, $f_4(s)$ and $f_5(s)$ be polynomials of degree
$3$, $4$ and $5$ with no zero in $s=0$.
Then we substitute
$$
  t\,=\,\frac{s^6}{f_5(s)},\,t\,=\,\frac{s^5}{f_4(s)}\,\mbox{ \quad and \quad } 
  t\,=\,\frac{s^4}{f_3(s)}
$$ 
into the Weierstra\ss\ equation $y^2+txy=x^3-t^5$ of the
universal family over $\oIgord{3}$, see Proposition \ref{explicit equations}.
In all cases this leads to a Weierstra\ss\ equation 
$$
  y^2\,=\,x^3\,+\,s^2x^2\,+\,s^{5}\,+\,r_4s^{4}\,+\,r_3s^{3}\,+\,r_2s^{2}\,+\,r_1s\,+\,r_0
$$
for certain $(r_4,r_3,r_2,r_1,r_0)\in \Aff_k^5$.
Depending on the degree of $\varphi$ these coefficients
satisfy the following conditions:
$$\begin{array}{lcll}
 \deg\varphi\,=\,6 &:& r_1 r_0\,\neq\,0 \\
 \deg\varphi\,=\,5 &:& r_1\,\neq\,0,& r_0\,=\,0\\
 \deg\varphi\,=\,4 &:& r_2\,\neq\,0, & r_1\,=\,r_0\,=\,0
\end{array}
$$
Note that the generic surfaces of each degree correspond to the extremal
rational surfaces of the cases 1C, 1D and 3C of \cite[\S3]{Lang94}.

It is remarkable that these rational elliptic surfaces appear in the
family of elliptic surfaces related to the semi-universal deformation of
the $E_8^2$-singularity in characteristic $3$, which is given by
$$
  y^2\,=\,x^3\,+\,(t^2+s)x^2\,+\,(q_1t+q_0)x\,+\,t^5\,+\,r_4t^4\,+\,
  r_3t^3\,+\,r_2t^2\,+\,r_1t\,+\,r_0.
$$

To obtain elliptic K3 surfaces with $3$-torsion section we have to 
take the Frobenius pullback of these surfaces.
Then the non-trivial $3$-torsion sections of the fibration
are explicitly given by
$$
(-(\tau ^5+r_4^{\frac{1}{3}}\tau ^4+r_3^{\frac{1}{3}}\tau
^3+r_2^{\frac{1}{3}} \tau ^2+r_1^{\frac{1}{3}} \tau
+r_0^{\frac{1}{3}}), 
\pm \tau ^3 (\tau ^5+r_4^{\frac{1}{3}}\tau ^4+r_3^{\frac{1}{3}}\tau
^3+r_2^{\frac{1}{3}} \tau ^2+r_1^{\frac{1}{3}} \tau
+r_0^{\frac{1}{3}}))
$$
(For $\deg \varphi=4$ one needs to modify slightly
because of the minimality of the equation.)

By Lemma \ref{determinant formula} and the preceding argument, 
the index of $\MW(Y)_{\rm free}(3)$
inside $\MW(X)_{\rm free}$ is related to 
the Artin invariant of $X$ for each case in the table. 
From this observation we obtain an upper bound for the Artin
invariant. 
On the other hand, 
since all the surfaces in the table can be realized inside the family
corresponding to the semi-universal deformation of the $E_8^2$-singularity as
noted above, 
the dimension of the surface having the given type of singular fibers
inside the moduli space is bounded from below.
This gives the Artin invariants for the cases  
$\deg\varphi = 4$ and $\deg \varphi =6$.

For the case $\deg \varphi =5$ we need a more precise analysis. 
Let $X$ be an elliptic K3 surface with $3$-torsion sections whose
singular fibers are of type ${\rm IV}_5, 5\times {\rm I}_3$. 
Then we have $\mu^2 = 3^{12-2\sigma_0(X)}$, where $\mu$ is
the index of $\MW (Y)_{\rm free}(3)$ inside $\MW(X)_{\rm free}$.
This implies $\sigma_0(X)\leq 6$.
On the other hand, these surfaces are realized inside
the semi-universal deformation of the $E_8^2$-singularity, which
yields $\sigma_0(X)\geq 5$.
Thus, we have to decide whether $\mu=1$ or $\mu=3$ holds true.
Assume $\mu=1$.
From $\MW(Y)_{\rm free}=\nMW (Y)=E_8$ 
we get $\MW(X)_{\rm free}=\nMW(X)=E_8(3)$. 
However, the $3$-torsion sections of this surface do not lie
in $\nMW(X)$,
which produces many free sections in $\MW (X)$ that do not lie in
$\nMW (X)$, a contradiction.
Thus, $\mu=3$ and we obtain $\sigma_0(X)=5$.
The other cases can be treated similarly using
Lemma \ref{tech lemma 3}.

Since we have found examples for all Artin invariants 
$\sigma_0\leq6$, Proposition \ref{complete families} tells
us that every (Shioda-)supersingular K3 surface with $\sigma_0\leq6$
possesses an elliptic fibration with $3$-torsion section.
\qed

\begin{Lemma}
  \label{tech lemma 3}
  With the notations as before,
  the index of $\nMW(Y)$ inside $\MW (Y)_{\rm free}$ divides
  the index of $\nMW(X)$ inside $\MW (X)_{\rm free}$.\qed
  \end{Lemma}

\section{Characteristic $2$}
\label{sec:p=2}

In this section we deal with elliptic K3 surfaces
with $2$-torsion section in characteristic $2$.
The classification in this case has much more 
subcases as for $p\geq3$ since
the fibration may have constant $j$-invariant, additive fibers 
may not be potentially supersingular
and potentially supersingular may have good reduction.

We start with a useful result, which directly follows from
\cite{Dolgachev; Keum 2001}:

\begin{Proposition}
 \label{useful}
 Let $X\to\PP^1$ be an elliptic K3 surface with $2$-torsion section in 
 characteristic $2$.
 Then $X=Y^{(2)}$ for some elliptic fibration $Y\to\PP^1$.
 Moreover, denote by $G$ the group of order $2$ that acts on $X$ via
 translating by the $2$-torsion point.
 Then $Y=X/G$ and there are two cases
 \begin{enumerate}
  \item $G$ has one or two fixed points and $Y$ is a K3 surface
  \item The fixed locus of $G$ is a connected 
    curve and $Y$ is a rational surface.
    In particular, $X$ is unirational in this case.
 \end{enumerate}
\end{Proposition}

\proof
Let us recall that
multiplication by $2$ on generic fibers of the fibration factors as
$Y\to Y^{(2)}=X \to X/G=Y$, cf. (\ref{multiplication by p}).

If $G$ has a finite number of fixed points then there are at most two
of them by \cite[Theorem 2.4]{Dolgachev; Keum 2001}.
If $G$ acted without fixed points, then $Y$ would be an Enriques surface,
which is absurd, cf. the proof of Theorem \ref{at most two fibers}.
If $G$ has one fixed point then $X/G$ is a K3 surface by 
\cite[Theorem 2.4]{Dolgachev; Keum 2001} and
\cite[Remark 2.6]{Dolgachev; Keum 2001}.
And if $G$ has two fixed points then $X/G$ is also a K3 surface by 
\cite[Theorem 2.4]{Dolgachev; Keum 2001}.

If $G$ has non-isolated fixed points then the fixed locus is a connected
curve by \cite[Corollary 3.6]{Dolgachev; Keum 2001} and the quotient 
$X/G$ is rational \cite[Theorem 3.7]{Dolgachev; Keum 2001}.
\qed\medskip

The classification of elliptic K3 surfaces with $2$-torsion in characteristic
$2$ is now as follows, where $h$ denotes the height of the formal
Brauer group as discussed in Section \ref{sec:brauer}.

\begin{Theorem}
  \label{p=2theorem}
  Let $X\to\PP^1$ be an elliptic K3 surface with $2$-torsion section in 
  characteristic $2$.

  If the fibration has constant $j$-invariant then the singular fibers are either
  \begin{enumerate}
    \item one additive fiber of type ${\rm I}^*_{12,6}$, and then $h\geq2$, or
    \item two additive fibers, both of type ${\rm I}^*_{4,2}$,
       and then $h=1$.
  \end{enumerate}

  If the fibration does not have constant $j$-invariant, then we have the following cases:
  \begin{enumerate}
    \item the fibration has precisely one additive
       fiber, which is potentially supersingular.
       In this case $h\geq2$ holds true.

    \item the fibration is semi-stable and there is precisely one 
       fiber with good and supersingular reduction.
       Moreover, $X$ is unirational and $h=\infty$.

    \item the fibration has precisely two fibers with
        additive reduction, both of which are potentially supersingular.
        In this case $h=1$ holds true.

    \item the fibration has precisely two fibers with
        additive reduction, one of which is potentially supersingular
        and the other one is potentially ordinary of type
        ${\rm I}^*_{4,2}$.
        In this case $h=1$ holds true.
  \end{enumerate}
\end{Theorem}

\proof
Let
$$
y^2+a_1(t)xy+a_3(t)y=x^3+a_2(t)x^2+a_4(t)x+a_6(t)
$$
be a global Weierstra\ss\ equation of the K3 surface, where the $a_i(t)$'s are 
polynomials of degree $\leq 2i$.
We denote by $\sigma_2$ the $2$-torsion section and let $G$ be the group of
order $2$ generated by translation by $\sigma_2$.

In order to have additive reduction at $t_0$ it is necessary that $a_1(t_0)=0$.
As $\deg a_1(t)\leq2$, it follows that there are at most two places of
additive reduction.
Moreover, from 
$j=a_1^{12}/\Delta$
we infer that for a place $t_0$ to have potentially supersingular 
reduction, again, $a_1(t_0)=0$ is necessary.

{\sc Case 1}: assume that $a_1(t)$ has a double zero. 
Then we get $h\geq2$ from \cite[Theorem 2.12]{Artin 1974}.

If the fibration does not have constant $j$-invariant then there has to be at least one place of
potentially supersingular reduction, which corresponds to the
double zero of $a_1(t)$.
Then this fiber has either additive reduction and we are in case (1)
or else this fiber has good supersingular reduction and we are in
case (2).
In this latter case $(\sigma_2\cdot\sigma_0)\geq1$,
translation by $\sigma_2$ fixes the whole supersingular fiber
and the quotient $X/G$ is rational by Proposition \ref{useful}.
In particular, $X$ is unirational and thus $h=\infty$.

If the fibration has constant $j$-invariant then the only singular
fibers can be of type ${\rm I}^*_{4+8d}$ for some $d\geq0$, which
have minimal discriminant $12d+12$ and Swan conductor $2+4d$
by \cite[Proposition 15.1]{Liedtke; Schroeer 2008}.
Since the minimal discriminants add up to $c_2(X)=24$ and there is at most
one additive fiber it has to be of type ${\rm I}^*_{12,6}$.

{\sc case 2}: assume that $a_1(t)$ has two distinct zeros.
From \cite[Theorem 2.12]{Artin 1974} we obtain $h=1$.
In particular, $X/G$ is a K3 surface and the $G$-action
has one or two fixed points by Proposition \ref{useful}.

In the case where $j(t)$ is constant it has to be a unit at both simple
zeros of $a_1(t)$, i.e., $\Delta$ has a zero of order $12$ at both places.
Depending on whether the discriminant is minimal, the reduction at such
a place is either good or of type ${\rm I}^*_{4,2}$ by
\cite[Proposition 15.1]{Liedtke; Schroeer 2008}.
Since the sum of the minimal discriminants is equal to $c_2(X)=24$
we must have two fibers of type ${\rm I}^*_{4,2}$.

We may thus assume that the fibration does not have constant
$j$-invariant.

First, assume that both places are potentially supersingular.
Then both places have additive reduction since the
$G$-action has two fixed points and would fix a supersingular
fiber completely by Proposition \ref{fixed locus of translation}.
This is case (3).

Now, assume that one of the zeros of $a_1(t)$ corresponds to a place
with potentially ordinary or potentially multiplicative reduction.
Not both zeros can belong to places of good or multiplicative reduction
since there is at least one potentially supersingular fiber.
Let $t_0$ be the place with potentially good or ordinary reduction.
By \cite[Section 15]{Liedtke; Schroeer 2008} the minimal discriminant 
at this place equals
$$
v(\Delta) = 12 + 12d - 2v_{t_0}(j) \geq 12
$$
and the reduction is of type ${\rm I}^*_{4+8d-2v(j)}$
As explained in loc. cit. such an additive and not potentially
supersingular fiber arises a quadratic twist from an 
elliptic fibration $X'\to\PP^1$
that has semi-stable reduction at the place corresponding to the
${\rm I}^*_{4+8d-2v(j)}$-fiber.
This quadratic twist may be arranged in 
such a way that the other fibers are not affected, which
implies $c_2(X')<c_2(X)$.
Since $j(X)=j(X')$ the fibration still has non-constant
$j$-invariant after twisting and thus $c_2(X')\neq0$.
In particular $c_2(X')=12$, i.e., $X'$ is a rational surface.
This implies that the minimal discriminant of the fiber of type
${\rm I}^*_{4+8d-2v(j)}$ equals $12$, i.e., $d=0$ and $v(j)=0$ and
we get a fiber of type ${\rm I}^*_{4,2}$ with potentially ordinary
reduction.
Also the potentially supersingular fiber must have additive
reduction or else the $G$-action would fix a supersingular fiber,
but we already now that $G$ fixes only two points.
\qed\medskip

\begin{Remark}
  Compared to characteristic $p\geq3$ the new, ``exotic'' 
  classes are fibrations with constant $j$-invariant, as well 
  as classes (2) and (4) in the case of non-constant $j$-invariant.
  We will classify them completely in Section \ref{sec:newp=2}.
  There, it will turn out that they are supersingular if only if the
  elliptic fibration arises as Frobenius pullback from a rational 
  surface, as predicted by Artin's Conjecture \ref{Artin Conjecture2}.
\end{Remark}

In characteristic $p\geq3$, an elliptic K3 surface with $p$-torsion
section that has precisely one fiber with
potentially supersingular reduction is supersingular, unirational 
and its elliptic fibration arises as Frobenius pullback from 
a rational elliptic surface.

The following examples have one additive and potentially supersingular
fiber, i.e., the height $h$ of the formal Brauer group is at least
$2$ by Theorem \ref{Brauer height}.
However, these surfaces are not supersingular and their elliptic
fibrations arise as Frobenius pullback from K3 surfaces, i.e.,
the alternative of Theorem \ref{ordinary characterisation}
does not hold in characteristic $2$.

\begin{Proposition}
 Let $X\to\PP^1$ be the elliptic K3 surface given by the
 Weierstra\ss\ equation  
 $$
  y^2\,+\,t^2xy\,+\,t^2y\,=\,x^3\,+\,(1+t)x^2\,+\,t\,.
 $$
 The elliptic fibration is a $4$-fold Frobenius pullback.
 More precisely,
 $$
 \begin{array}{ccccc}
  & j-\mbox{invariant} & \mbox{singular fibers} & \mbox{type} & \mbox{ height of $\Brauer$}\\
  \hline
  X & t^{16} & {\rm II}_6, {\rm I}_{16} & \mbox{K3} & 2\\
  X^{(1/2)} & t^8 & {\rm I}^*_{4,6}, {\rm I}_8 & \mbox{K3} & 2\\
  X^{(1/4)} & t^4 & {\rm I}^*_{8,6}, {\rm I}_4 & \mbox{K3} & 2\\
  X^{(1/8)} & t^2 & {\rm I}^*_{10,6}, {\rm I}_2 & \mbox{K3} & 2\\
  X^{(1/16)} & t & {\rm I}^*_{11,6}, {\rm I}_1 & \mbox{K3} & 2\\
 \end{array}$$
 The elliptic fibrations of
 $X$, $X^{(1/2)}$, $X^{(1/4)}$ and $X^{(1/8)}$
 possess $2$-torsion sections and
 arise as Frobenius pullbacks from K3 surfaces.
\end{Proposition}

\proof
The computation of the singular fibers is straight forward
and left to the reader.
Moreover, all surfaces are K3 surfaces and since they
are related by Frobenius pullbacks the heights of
their formal Brauer groups coincide.
Thus, it suffices to compute the formal Brauer group of one surface
and we take the one of the statement of the proposition.
Making a coordinate change to achieve $a_2=0$ in the
Weierstra\ss\ equation we can apply
\cite[Theorem (2.12)]{Artin 1974} and obtain $h=2$.
\qed\medskip

These surfaces belong to class (1) with non-constant
$j$-invariant of Theorem \ref{p=2theorem}.
We shall see further examples with $h=2$ and iso-trivial
fibrations in the next section.

\section{The exotic classes in characteristic $2$}
\label{sec:newp=2}

This sections deals with the classes of
Theorem \ref{p=2theorem} that do no exist
for $p\geq3$.

\subsection*{Fibrations with constant $\mathbf{j}$-invariant}
This class coincides with the Kummer surfaces 
studied by Shioda in \cite{Shioda 1974}:

\begin{Proposition}
 Every elliptically fibered K3 surface with constant $j$-invariant
 and $2$-torsion
 section in characteristic $2$ arises as minimal
 desingularization of
 \begin{equation}
   \label{isotriv}
   (E_1\times E_2)/G \,\to\, E_2/G\,\iso\,\PP^1\,,
 \end{equation}
 where $E_1$ is an ordinary and $E_2$ is an arbitrary 
 elliptic curve, and $G\iso\ZZ/2\ZZ$ acts via the
 sign involution on each factor.

 Conversely, for any two elliptic curves $E_1$, $E_2$,
 where $E_1$ is ordinary, a minimal desingularization 
 of (\ref{isotriv}) yields an 
 elliptic K3 surface with constant $j$-invariant and
 $2$-torsion section.
 More precisely,
 $$\begin{array}{lccc}
      E_2 & \mbox{singular fibers} & \rho & h\\
      \hline
      \mbox{ordinary} & 2\times{\rm I}_{4,2}^* & 18\leq\rho\leq20 & 1\\
      \mbox{supersingular} & {\rm I}_{12,6}^* & 18 & 2
   \end{array}$$ 
 In particular, these surfaces cannot be supersingular, 
 and $h=2$ is possible.
\end{Proposition}

\proof
Since the generic fiber is ordinary, such a surface is a quadratic 
twist of a trivial fibration.
Thus, $X$ arises via $(E_1\times C)/G\to C/G$,
where $\varphi:C\to\PP^1$ 
is an Artin--Schreier morphism of degree $2$.
The group $G=\ZZ/2\ZZ$ acts via the sign involution on $E_1$ and
via the Galois action on $C$.
From Theorem \ref{p=2theorem} we know that the fibration $X\to\PP^1$ 
has either one fiber of type ${\rm I}^*_{12,6}$ or two fibers
of type ${\rm I}^*_{4,2}$.
From \cite[Section 15]{Liedtke; Schroeer 2008} it then follows
that $\varphi$ is ramified in one point with four non-trivial
higher ramification groups (the ${\rm I}_{12,6}^*$-case)
or in two points with two non-trivial higher ramification groups
(the $2\times{\rm I}_{4,2}^*$-case).
In both cases $C$ is an elliptic curve, and the Galois action
coincides with the sign involution.
In case, $\varphi$ is ramified in one point, its $p$-rank is
trivial \cite[Corollary 1.8]{Crew 1984}, and thus $C$ is
supersingular.
Similarly, if $\varphi$ is ramified in two points then $C$ is
ordinary.

Conversely, it is easy to see that this construction yields
elliptic K3 surfaces with $2$-torsion section.

The rank $\rho$ of the N\'eron--Severi group has been determined
in \cite{Shioda 1974}.
We set $A=E_1\times E_2$, where $E_1$ is an ordinary elliptic curve.
Then the height of $\Brauer(A)$ is $1$ or $2$ depending on whether
$E_2$ is ordinary or supersingular 
\cite[Lemma 6.2]{Geer; Katsura 2003}.
Since $A/G$ has only rational singularities  
\cite{Shioda 1974}
we can conclude as in the proof of
\cite[Theorem 6.1]{Geer; Katsura 2003} that the 
formal Brauer groups of $A/G$ and $X$ are isomorphic.
Since $A\to A/G$ is an Artin-Schreier covering of degree $2$,
there is a non-trivial trace map, and as in the proof of
\cite[Theorem 6.1]{Geer; Katsura 2003} we conclude that
the formal Brauer groups of $A$ and $A/G$ are isomorphic.
\qed

\subsection*{Semi-stable fibrations}
Class (2) with non-constant $j$-invariant in
Theorem \ref{p=2theorem} is closely related to rational
elliptic surfaces.
These surfaces are unirational and supersingular.

\begin{Proposition}
  \label{p=2semistable}
  Let $X\to\PP^1$ be an elliptic K3 surface with $2$-torsion
  section in characteristic $2$ whose fibration is semi-stable.
  Then $X\to\PP^1$ arises as Frobenius pullback from a rational
  elliptic surface $Y\to\PP^1$ with semistable fibration.
  
  Conversely, if $Y\to\PP^1$ is a rational elliptic surface 
  with semistable fibration, then its Frobenius pullback yields
  an elliptic K3 surface with $2$-torsion section.
\end{Proposition}

\prf
We have seen in the proof of Theorem \ref{p=2theorem} 
that $Y\to\PP^1$ is rational.
Moreover, the elliptic fibration on $Y$ must be semi-stable
because the one on $X$ is.
We leave the converse to the reader.
\qed

\begin{Remark}
  In \cite[Section 4]{Ito Math Nachr}, an $8$-dimensional family 
  of semistable rational elliptic surfaces related to the  
  deformation of an $E_8^4$-singularity is constructed.
  Via Frobenius pullback we obtain an 
  $8$-dimensional family of semistable elliptic K3 surfaces
  with Artin invariants $1\leq\sigma_0\leq9$, see
  \cite[Theorem 5.2]{Ito Math Nachr}.
\end{Remark}

\subsection*{Additive and potentially ordinary fibers}
Also, Class (4) with non-constant $j$-invariant in
Theorem \ref{p=2theorem} is closely related to rational
elliptic surfaces.
However, being ordinary, these surfaces are neither 
unirational nor supersingular.

In order to state the result, let us introduce
the following notation:
For a point $Q\in\PP^1$ denote by
$\psi_Q:\PP^1\to\PP^1$ the Artin--Schreier morphism
of degree $2$ that is branched over $Q$.

\begin{Proposition}
  Let $X\to\PP^1$ be an elliptic K3 surface
  with non-constant $j$-invariant and 
  $2$-torsion section in characteristic $2$ that 
  possesses a potentially ordinary fiber of 
  type ${\rm I}^*_{4,2}$, say, at $Q\in\PP^1$.
  Then there exists a rational elliptic surface 
  $X'\to\PP^1$ with $2$-torsion section and 
  good ordinary reduction at $Q$
  such that $X$ arises as quadratic twist from $X'$
  via $\psi_Q$.

  Conversely, if $X'\to\PP^1$ is a rational elliptic surface
  with $2$-torsion section and with good ordinary reduction
  at $Q\in\PP^1$ then the quadratic twist of
  $X'$ with respect to $\psi_Q$ yields
  an elliptic K3 surface with $2$-torsion section
  and a potentially ordinary fiber of type ${\rm I}^*_{4,2}$
  above $Q$.
\end{Proposition}

\proof
From \cite[Section 15]{Liedtke; Schroeer 2008} we see that
the ${\rm I}^*_{4,2}$-fiber arises from an elliptic
fibration $X'\to\PP^1$ as quadratic twist
$\psi:C\to\PP^1$, which is totally ramified at $Q$.
If $Q$ is the only branch point, which we can and will assume,
then $X'$ has the same singular fibers as $X$ but has good reduction at
$Q$.
In particular, $c_2(X')<c_2(X)$, and so 
$X'\to\PP^1$ is a rational elliptic surface.
Since the reduction type of $X\to\PP^1$ at $Q$ has
Swan conductor $\delta=2$, we conclude that $\psi$ has two
non-trivial higher ramification groups, i.e., $\psi=\psi_Q$.
We leave the converse to the reader.
\qed

\begin{Remark}
 Note that $X$ and $X'$ have the same numbers and types
 of singular fibers (including Swan conductors) except
 for the ${\rm I}^*_{4,2}$-fiber at $Q$ which is induced 
 on $X$ by the quadratic twist.
\end{Remark}

\begin{Example}
To illustrate this case with an example, consider the 
universal elliptic curve ${\cal E}\to\Igord{4}$.
Then ${\cal E}^{(2)}\to\Igord{4}$ corresponds to a 
rational elliptic surface with $2$-torsion section.
Twisting with respect to $\psi_Q:\PP^1\to\PP^1$, where 
$Q\in\Igord{4}\subset\PP^1$ corresponds to the
ordinary $j$-value $j=1$, we obtain
$$
y^2\,+\,xy\,=\,x^3\,+\,\frac{1}{t+1}x^2\,+\,t^2\,.
$$
This is an elliptic K3 surface with $2$-torsion section,
having singular fibers of type ${\rm I}_2$, 
${\rm III}^*_1$ (potentially supersingular) 
and ${\rm I}^*_{4,2}$ (potentially ordinary).
\end{Example}

\section{sections of order $4$ and $8$}
\label{sec:p=4}

In this section we classify elliptic K3 surfaces with 
$8$- and $4$-torsion sections in characteristic $2$, using
again the Igusa curves.
They turn out to belong to classes (1) and (3) with non-constant
$j$-invariant of Theorem \ref{p=2theorem}.
It also turns out that Theorem \ref{supersingular characterisation}
and Theorem \ref{ordinary characterisation} hold for them.
Thus, these surfaces behave like the ones in characteristic
$p\geq3$.

\subsection*{$\mathbf{8}$-torsion sections}
The following result is proved as Theorem \ref{p=7}, which is why we
leave it to the reader:

\begin{Theorem}
  \label{p=8}
  There exists only one elliptic K3 surface 
  $X\to\PP^1$ with $8$-torsion section in characteristic $2$ 
  up to isomorphism.
  It has the following invariants:
  $$\begin{array}{lccc}
     \mbox{singular fibers} & \sigma_0 & \nMW(X) & \MW(X) \\
     \hline
      {\rm I}^*_{1,1}, 2\times{\rm I}_8 & 1 & A_1(2) & A_1^*(2)\oplus(\ZZ/8\ZZ) 
    \end{array}
  $$
  The Weierstra\ss\ equation is given by the following:
  $$
     y^2\,+\,t^2xy\,=\,x^3\,+\,x\,+\,t^4.
  $$
  In particular, it is the unique supersingular K3 surface with 
  Artin invariant $\sigma_0=1$.
\end{Theorem}

\begin{Remark}
  Having Artin invariant $\sigma_0=1$, it is a
  generalized Kummer surface \cite{Schroeer 2007}.
  An explicit Weierstra\ss\ equation
  is given in Proposition \ref{explicit equations},
  but we note that
  uniqueness and an equation have already been
  obtained in \cite[Examples 2.4]{Schweizer 2005}.
\end{Remark}

\subsection*{$\mathbf{4}$-torsion sections}

As before, we denote by $O\in\oIgord{4}$ the unique supersingular point.
Since our proof works as for $p=3$ or $p=7$, we leave it to
the reader and only state the result:

\begin{Theorem}
  In characteristic $2$, the classifying morphism $\varphi$ 
  for an elliptic K3 surface with
  $4$-torsion section is finite of degree
  $2\leq\deg\varphi\leq4$. 
  More precisely,
  \begin{enumerate}
   \item $\deg\varphi=2$, $\varphi$ is separable and $\varphi^{-1}(O)$ consists of
    two points, or
   \item $\deg\varphi=3$ and $\varphi^{-1}(O)$ consists of one point
   or two points, or
   \item $\deg\varphi=4$ and $\varphi^{-1}(O)$ consists of one point or two points
      with ramification index $e=2$ (wildly ramified).
  \end{enumerate}
  Conversely, if $\varphi$ is as above then the associated elliptic fibration
  with $4$-torsion section is a K3 surface.
\end{Theorem}

More precisely, depending on the branch points we obtain the following table,
where $X=Y^{(4)}$ and the type of $Y$ is tabled in the last column
$$\begin{array}{clllccc}
     \deg\varphi&\mbox{ \quad} &
     \multicolumn{2}{c}{\mbox{singular fibers}} & \dim & h & 
     Y \\
     \hline
      2 & \multicolumn{4}{l}{\varphi \mbox { separable:}}\\
        & & 2\times{\rm I}_{1,1}^* & 2\times {\rm I}_4 & 2 & 1 & \mbox{K3} \\
        & & 2\times{\rm I}_{1,1}^* & {\rm I}_8 & 1 & 1 & \mbox{K3} \\
     \hline
      3 & & {\rm I}_{1,1}^*, {\rm III}_1 & 
        3\times{\rm I}_4  
        & 3 & 1 & \mbox{K3} \\
        & & {\rm I}_{1,1}^*, {\rm III}_1 & 
        {\rm I}_8, {\rm I}_4
        & 2 & 1 & \mbox{K3} \\
        & & {\rm I}_{1,1}^*, {\rm III}_1 & 
        {\rm I}_{12}
        & 1 & 1 & \mbox{K3} \\
        & & {\rm I}_{3,3}^* &
        3\times{\rm I}_4 
        & 2 & \infty & \mbox{rational}\\
        & & {\rm I}_{3,3}^* &
        {\rm I}_8, {\rm I}_4
        & 1 & \infty & \mbox{rational}\\
        & & {\rm I}_{3,3}^* &
        {\rm I}_{12}
        & 0 & \infty & \mbox{rational}\\
     \hline
      4 & \multicolumn{4}{l}{\varphi \mbox { separable:}}\\
        & & 2\times {\rm III}_1 & 
        4\times {\rm I}_4
        & 4 & 1 & \mbox{K3} \\
        & & 2\times {\rm III}_1 & 
        {\rm I}_8, 2\times {\rm I}_4
        & 3 & 1 & \mbox{K3} \\
        & & 2\times {\rm III}_1 & 
        {\rm I}_{12}, {\rm I}_4
        & 2 & 1 & \mbox{K3} \\
        & & {\rm I}_{0,2}^* & 
        4\times{\rm I}_4 
        & 3 & \infty & \mbox{rational}\\
        & & {\rm I}_{0,2}^* & 
        {\rm I}_8, 2\times{\rm I}_4
        & 2 & \infty & \mbox{rational}\\
        & & {\rm I}_{0,2}^* & 
        {\rm I}_{12}, {\rm I}_4 
        & 1 & \infty & \mbox{rational}\\
        & \multicolumn{4}{l}{\varphi \mbox { inseparable but not purely inseparable:}}\\
        & & 2\times {\rm III}_1 & 
        2\times{\rm I}_{8}
        & 2 & 1 & \mbox{K3}\\
        & & 2\times {\rm III}_1 & 
        {\rm I}_{16}
        & 1 & 1 & \mbox{K3}\\
        & & {\rm I}_{1,1}^* & 2\times{\rm I}_{8} 
        & 0 & \infty & \mbox{rational}\\
        & \multicolumn{4}{l}{\varphi \mbox { purely inseparable:}}\\
        & & {\rm I}_{1,1}^* & {\rm I}_{16} & 0 & \infty & \mbox{rational}
 \end{array}$$


\begin{Remark}
  There are two unique surfaces in this list:
  \begin{itemize}
    \item[-] If $\deg\varphi=4$,
     $\varphi$ is inseparable but not purely inseparable and
     $\varphi$ totally ramified over $O$,
     we obtain the unique surface with $8$-torsion section.
    \item[-]  
     If $\deg\varphi=4$ and $\varphi$ is purely inseparable,
     the resulting elliptic K3 surface is extremal.
     Such surfaces in characteristic $p=2,3$ have been studied 
     and classified in \cite{Ito 2002}.
     In fact, our surface appears in Table 1 of loc. cit.
  \end{itemize}
\end{Remark}

Similar to characteristic $3$, the generic supersingular surfaces can be related to
deformations of singularities.
Namely, 
$$
 y^2\,+\,txy\,=\,x^3\,+\,t^5\,+\,r_4t^4\,+\,r_3t^3\,+\,r_2t^2
$$
defines a $3$-dimensional family of rational elliptic
surfaces.
This family arises as subfamily of 
the semi-universal deformation of a $E_8^4$-singularity.
Then all $Y^{(4)}_\lambda\to\PP^1$ are elliptic
K3 surfaces with $4$-torsion sections.
We leave the following result, whose proof is analogous
to the one of Theorem \ref{char. 3 table} to the reader.

\begin{Theorem}
  Every (Shioda-)supersingular K3 surface with Artin invariant
  $\sigma_0\leq4$ in characteristic $2$ possesses an elliptic
  fibration with $4$-torsion section.

  The complete list of these surfaces is given by the 
  following table
$$\begin{array}{clllcccc}
     \deg\varphi&\mbox{ \quad} &
     \multicolumn{2}{c}{\mbox{singular fibers}} & \dim & \sigma_0 &
     \nMW(X) & \MW(X) \\
     \hline
      3 & & {\rm I}_{3,3}^* &
        3\times{\rm I}_4 
        & 2 & 3 & D_4(2) & D_4^*(2)\oplus {\ZZ/4\ZZ} \\
        & & {\rm I}_{3,3}^* &
        {\rm I}_8, {\rm I}_4
        & 1 & 2 & A_3 & A_3^*\oplus {\ZZ/4\ZZ} \\
        & & {\rm I}_{3,3}^* &
        {\rm I}_{12}
        & 0 & 1 & A_2 & A_2^*\oplus {\ZZ/4\ZZ} \\
     \hline
      4 & \multicolumn{5}{l}{\varphi \mbox { separable:}}\\
        & & {\rm I}_{0,2}^* & 
        4\times{\rm I}_4 
        & 3 & 4 & D_4(2) & D_4^*(2)\oplus {\ZZ/4\ZZ}\\
        & & {\rm I}_{0,2}^* & 
        {\rm I}_8, 2\times{\rm I}_4
        & 2 & 3 & A_3(2) & A_3^*(2)\oplus {\ZZ/4\ZZ}\\
        & & {\rm I}_{0,2}^* & 
        {\rm I}_{12}, {\rm I}_4 
        & 1 & 2 & A_2(2) & A_2^*(2)\oplus {\ZZ/4\ZZ}\\
        & \multicolumn{5}{l}{\varphi \mbox { inseparable but not purely inseparable:}}\\
        & & {\rm I}_{1,1}^* & 2\times{\rm I}_{8} 
        & 0 & 1 & A_1(2) & A_1^*(2)\oplus {\ZZ/8\ZZ}\\
        & \multicolumn{5}{l}{\varphi \mbox { purely inseparable:}}\\
        & & {\rm I}_{1,1}^* & {\rm I}_{16} & 0 & 1 & \{ 0\} & {\ZZ/4\ZZ}
 \end{array}$$

\end{Theorem}

Moreover, from the table above we see that the implications
(4) $\Rightarrow$ (5) $\Rightarrow$ (1) of 
Theorem \ref{supersingular characterisation}
hold for these surfaces.
Also, we see that these surfaces can only be ordinary or supersingular.
Thus,

\begin{Theorem}
  Theorem \ref{supersingular characterisation} and 
  Theorem \ref{ordinary characterisation}
  hold for elliptic K3 surfaces with $4$-torsion section in characteristic $2$.
\end{Theorem}

\end{document}